%% file: main.tex
\documentclass[twoside,11pt]{article}

\input{packages}
\input{commands}

\jmlrheading{24}{2025}{1-\pageref{LastPage}}{7/24}{-}{}{Florian Heinrichs, Patrick Bastian, Holger Dette}

\ShortHeadings{Sequential Outlier detection in non-stationary time series}{F. Heinrichs, P. Bastian, H. Dette}
\firstpageno{1}

\begin{document}
	
\title{Sequential Outlier Detection in Non-Stationary Time Series}

\author{\name Florian Heinrichs \email f.heinrichs@fh-aachen.de \\
	\addr FH Aachen\\
	Heinrich-Mußmann-Straße 1\\
	52428 Jülich, Germany 
    \AND
    \name Patrick Bastian \email patrick.bastian@rub.de \\
	\addr Ruhr-Universität Bochum\\
	Universitätsstraße 150\\
	44801 Bochum, Germany 
	\AND
	\name Holger Dette \email holger.dette@rub.de \\
	\addr Ruhr-Universität Bochum\\
	Universitätsstraße 150\\
	44801 Bochum, Germany
}

\editor{-}

\maketitle

\begin{abstract}%
	A novel method for sequential outlier detection in non-stationary time series is proposed. The method tests the null hypothesis of ``no outlier'' at each time point, addressing the multiple testing problem by bounding the error probability of successive tests, using extreme value theory. The asymptotic properties of the test statistic are studied under the null hypothesis and alternative. The finite sample properties of the new detection scheme are investigated by means of a simulation study, and the method is compared with alternative procedures  which have recently been proposed in the statistics and machine learning literature.
\end{abstract}

\begin{keywords}
	Outlier detection, Non-stationary time series, Local linear regression, Extreme value theory
\end{keywords}

\begin{msc_class}
    62M10, 62G10, 62L10
\end{msc_class}

\maketitle 

\input{sec_intro}

\input{sec_related_work}

\input{sec_methodology}

\input{sec_empirical_results}

\section{Proofs}

\input{app_main_thm}

\input{app_test_properties}

\input{app_jackknife_estimator}

\bibliography{bibliography}
	
\appendix

\newgeometry{left=0.7in, right=0.7in, top=1in, bottom=1in}

\input{app_empirical_results}

\end{document}

%% file: packages.tex
\usepackage{blindtext}

\usepackage[preprint]{jmlr2e}
\usepackage[a4paper]{geometry}

\usepackage[utf8]{inputenc}
\usepackage[T1]{fontenc}
\usepackage{lmodern}
\usepackage{amsmath}
\usepackage{paralist}
\usepackage{float}
\usepackage{accents}
\usepackage{color}
\usepackage{bm}
\usepackage{lscape}
\usepackage{subfig}
\usepackage{longtable}
\usepackage{multibib}
\usepackage{booktabs}
\newcites{suppl}{References} 
\usepackage{tikz}
\usetikzlibrary{positioning}
\usepackage{dsfont}
\usepackage{pgfplotstable}
\usepackage{hyperref}
\usepackage{lastpage}
\usepackage{setspace}

%% file: commands.tex
\newcommand{\eps}{{\varepsilon}}
\renewcommand{\phi}{\varphi}
\newcommand{\R}{\mathbb{R}}
\newcommand{\Z}{\mathbb{Z}}
\newcommand{\N}{\mathbb{N}}
\newcommand{\pr}{\mathbb{P}}       
\newcommand{\ex}{\mathbb{E}}       
 
\newcommand{\cov}{\textnormal{Cov}}

\newcommand{\Nc}{\mathcal{N}}
\newcommand{\Uc}{\mathcal{U}}

\newcommand{\Gc}{\mathcal{G}}
\newcommand{\Oc}{\mathcal{O}}

\newcommand{\diff}{{\,\mathrm{d}}}

\newcommand{\convw}{\rightsquigarrow}

\DeclareMathOperator*{\argmin}{argmin}

\newtheorem{assumption}[theorem]{Assumption}

%% file: sec_intro.tex
\section{Introduction} \label{sec:intro}

An important source of time-dependent data that is becoming increasingly widespread are sensors. However, when time series are recorded by sensors, a variety of errors can occur, for example incorrect or missing values due to measurement and transmission errors, sensor calibration issues, sensor drift, malfunctioning hardware, interference from environmental factors or battery failure. Many statistical methods are not robust and strongly influenced by such erroneous data, and it  is important to detect anomalies in time series automatically and reliably so that a downward statistical analysis leads to unbiased results. 

In some applications, the anomalies themselves are of interest because they allow conclusions about the underlying system. A typical example is an unexpected change of  the power consumption of a machine. In this case, the anomaly might not be caused by the sensor, but by a machine malfunction that requires human intervention.

The detection of \textit{outliers} or \textit{anomalies} in time series has a long history since the seminal paper by \cite{fox1972} and the early work by \cite{burman1988}. Nowadays, numerous methods for the detection of outliers in univariate, multivariate and functional data exist. 
Informally, outliers are often described as ``[...] an observation which deviates so much from the other observations as to arouse suspicions that it was generated by a different mechanism'' \citep{hawkins1980}. More formally, outliers are commonly modeled as \textit{additive outliers}, where a contaminated time series 
\begin{align} \label{dette1}Y_i = X_i + \delta \xi_i
\end{align}
is observed, for $i\in \N$. This type of outlier was considered, for instance, in  \cite{abraham1979}, \cite{Denby1979} and \cite{Dhar1991}. In model \eqref{dette1}, 
$(X_i)_{i\in\N}$
denotes the time series of interest, $\delta$ is a level of contamination, and $\xi_i$ indicates the existence of an outlier, where $\xi_i= 1$ if an outlier occurs at time $i$ and $\xi_i=0$ otherwise. In this case, we are interested in the sequential hypotheses
\begin{align} \label{dette2}
H_0^{(i)}: Y_i \sim F_{X_i} \quad \text{vs.} \quad H_1^{(i)}: Y_i \nsim F_{X_i}, 
\end{align}
for $i\in \N$, where $F_{X_i}$ denotes the cumulative distribution function of $X_i$. A common assumption in this model is that the time series $(X_i)_{i\in\N}$ is stationary. We relax this assumption by allowing the mean to change over time. 

The time series of interest, $(X_i)_{i\in\N}$, can be decomposed into a deterministic term and a random error, which  yields the location model $X_i = \mu_i + \eps_i$, where $\mu$ denotes the (unknown) mean function and $(\eps_i)_{i\in\Z}$ is a sequence of centered errors. Now, if we have a suitable estimator $\hat{\mu}$ of $\mu$, we can estimate the residuals $\hat{\eps}_i = Y_i - \hat{\mu}_i$. If the residual is larger than a critical value $c$, for some time point $i$, we reject the null hypothesis $H_0^{(i)}$ in \eqref{dette2} and the corresponding observation can be considered an outlier. The asymptotic properties of the resulting sequential testing procedure depend crucially on the critical value $c$.

A naive approach is to bound the probability of falsely rejecting $H_0^{(i)}$ by Chebychev's inequality. The resulting critical value is only a weak bound and expected to yield a conservative test. Despite this fact, the method is commonly used, especially in the field of machine learning, where outlier detection for time series has gained popularity in recent years. Different types of neural networks, such as Long Short-Term Memory networks (LSTMs), Convolutional Neural Networks (CNNs), Graph Neural Networks (GNNs), and Spiking Neural Networks (SNNs) have been used to estimate the mean function $\mu$ and reject $H_0^{(i)}$ for large residuals, based on Chebychev's inequality \citep{malhotra2015, munir2018, zhao2020, buchhorn2023, cherdo2023}.

With this approach, the sequential testing procedure not only loses power because the choice of the  critical value yields a conservative decision, but also because the test levels are controlled for the tests individually and dependence is not taken into account. This leads to an accumulation of error probabilities, and either many false rejections of the null hypothesis, or if the test levels are corrected, to an overly conservative test.

Traditionally, for ``short'' time series, this has not been a major problem, but with the spread of sensors, that record data with a high frequency, the accumulation of error probabilities becomes a serious issue. For example, assume we have a sensor that measures $100$ samples per second for $12$ hours per day, which yields $4\,320\,000$ observations daily. If we want to bound the probability of a false positive decision to 5\% per day, we would have to adjust the level of each individual test accordingly to $1.16^{-6}\%$, making the test useless. Clearly, most sensors do not record independent observations and this dependence must be taken into account. In this paper we propose a sequential testing procedure based on extreme value theory, which takes advantage of the dependency, and is therefore less conservative.

%% file: sec_related_work.tex
\noindent
{\bf Related Work}  The literature on outlier detection is abundant and can be roughly grouped into two categories: methodology that controls the level of each individual test and methodology that controls the level of multiple subsequent tests.

A notable example of the former is the method proposed by \cite{campulova2018}, which considered the location-scale model $X_i = \mu_i + \sigma_i \eps_i$, where $\mu_i$ and $\sigma(i)$ denote the (unknown) mean and variance functions and $(\eps_i)_{i\in\N}$ is a sequence of i.i.d. errors. The authors  proposed to use kernel regression for the estimation of $\mu$, detect possible change points in $\sigma$, and define critical values through Chebyshev's inequality. As they used a two-sided kernel to estimate $\mu$ and a retrospective change point detection method, their approach cannot be used for the sequential detection of change points. For a recent review on outlier detection in time series, see \cite{blazquez2021}.

In the field of machine learning, it was first proposed to estimate $\mu$ through LSTMs \citep{malhotra2015}. Subsequently, other neural networks were used as well, such as CNNs \citep{munir2018}, GNNs \citep{zhao2020, buchhorn2023} and SNNs \citep{cherdo2023}. 

Methods that fall into the second category and control the level of multiple tests, are generally based on extreme value theory (EVT). A common assumption in EVT is that observations are independent and identically distributed (i.i.d.). Even though there are relaxations of this assumption in the EVT literature, the i.i.d. assumption is still common when EVT is applied to outlier detection. In an early work, \cite{roberts1999} proposed to use EVT to detect outliers in a Gaussian mixture model. Subsequently, \cite{dupuis2004} fit a generalized extreme value distribution to a given data set and define observations as outliers that are not well described by the model. A different line of work studies additive outliers for a specific type of non-stationary time series. \cite{vogelsang1999} used the $t$-statistic to identify single outliers in an $ARIMA(p, 1, q)$ model. Later, \cite{perron2003} expanded this approach to enable the detection of multiple outliers. Building on this, \cite{burridge2006} introduced an EVT-based framework, where observations were flagged as outliers if the gap between consecutive order statistics was large. Further extending these methods, \cite{astill2013} proposed a bootstrap procedure for $ARIMA(p, d, q)$ processes, that is robust in terms of the order of integration ($d=0, 1$). Finally, \cite{bhattacharya2019} proposed a robust version of the Hill estimator to estimate the (positive) tail index of the data generating distribution, and classify observations as outliers, if they are unlikely for the estimated tail index. This work was later generalized to arbitrary real-valued tail indices \citep{bhattacharya2023}.

A notable exception from the i.i.d. assumption, is the work by \cite{holevsovsky2018}. It is generally developed under the above-mentioned location-scale model with i.i.d. errors, yet an application of EVT to stationary time series was briefly discussed. Note however, that the approach is based on a two-sided kernel estimate of $\mu$ and retrospective change point detection, so that the resulting outlier detection method is not suitable for sequential outlier detection.

%% file: sec_methodology.tex
\section{Methodology} \label{sec:method}

\subsection{The General Testing Problem} \label{sec:general}

For a time series $(X_i)_{i\in\N}$, we are interested in detecting outliers that deviate from the ``normal'' behavior of the time series. To be precise, we consider the model
\begin{equation} \label{eq:model}
    X_i = \mu\big(\tfrac{i}{n}\big) + \eps_{i} + c_i \quad i\in\N,    
\end{equation}
where $\mu:\R\to\R$ denotes an unknown mean function, $n\in\N$ is the resolution with which $X_i$ is measured, $(\eps_i)_{i\in\Z}$ is a sequence of (possibly dependent) centered random variables and $c_i$ is a sequence of variables indicating a normal observation, for $c_i=0$, or an outlier, for $c_i\neq 0$. In this model, the hypotheses from \eqref{dette2} correspond to 
\begin{equation*}
	H_0^{(i)}: X_i \sim F_{\eps_i}\big(\cdot - \mu(\tfrac{i}{n})\big) \quad vs. \quad H_1^{(i)}: X_i \nsim F_{\eps_i}\big(\cdot - \mu(\tfrac{i}{n})\big).
\end{equation*}
Since we are not interested in arbitrary small values of $c_i$, we assume in the following that $c_i$ takes on sufficiently large values under the alternative, as specified in Corollary \ref{thm:test_properties} (iii).

When sequentially testing for various hypotheses, we run into the multiple testing problem and the probabilities of falsely rejecting some null hypothesis accumulate. If we do not adapt the level of each single test, we falsely classify normal observations as outliers. Conversely, if we adapt the level, the power of the tests vanish and we miss outliers. In the following, we derive a decision rule that considers sequential dependencies of the error process $(\eps_i)_{i\in\Z}$ and allows balancing type I and type II errors.

\subsection{Mathematical Preliminaries} \label{sec:prelim}

Throughout, we assume $\mu$ to be sufficiently smooth, which is further specified in Assumption \ref{assump:mu}. In order to estimate the unknown mean function, we use a one-sided version of local linear regression. More specifically, define the local linear estimator of $\mu$ and its derivative, as

\begin{align*}
	\Big(\hat{\mu}_{h_n}(t), \widehat{\mu'}_{h_n}(t)\Big) 
	= \argmin_{b_0, b_1\in \R} \sum_{i\in\N} \Big(X_{i} - b_0 - b_1 \big(\tfrac{i}{n}-t\big) \Big)^2 K\big(\tfrac{i-nt}{nh_n}\big) .
\end{align*} 
for some kernel $K$ and bandwidth $h_n\searrow 0$. In the following, we assume that $K$ has the support $[-1, 0]$, which means that $\hat{\mu}_{h_n}(t)$ only depends on observations of time points $i\le nt$ and we can calculate the estimator without knowing future values of $X_i$ for $i > nt$. To account for the bias of this estimator, we employ the Jackknife bias reduction technique proposed by \cite{schucany1977} and define
\[ \tilde{\mu}_{n}(t) = 2 \hat{\mu}_{h_n/\sqrt{2}}(t) - \hat{\mu}_{h_n}(t). \]

Given the estimator $\tilde{\mu}_{n}$, we can calculate the residuals $\hat{\eps}_i = X_i - \tilde{\mu}_{n}(\tfrac{i}{n})$, and reject the null hypothesis of no point outlier at time $i$ for large values of $|\hat{\eps}_i|$. In Section \ref{sec:point}, we derive critical values based on the asymptotic behavior of $\tilde{\mu}_{n}$, that allow the construction of a test with some pre-defined asymptotic $\alpha$-level for $H_0$ that is consistent against a broad class of alternatives.

Before we state regularity conditions to derive the required results, recall some basic notations from EVT. The Generalized Extreme Value (GEV) distribution is defined by 
\[ G_{\gamma_0, \mu_0, \sigma_0}(x) = \exp\bigg\{ - \bigg( 1 + \gamma_0 \frac{x-\mu_0}{\sigma_0} \bigg)^{-\frac{1}{\gamma_0}} \bigg\},\quad  1 + \gamma_0 x > 0 \]
with parameters $\gamma_0$ (shape), $\mu_0$ (location) and $\sigma_0$ (scale). Note that $G_{0, \mu_0, \sigma_0}$ is interpreted as the limit 
\[
G_{0, \mu_0, \sigma_0}(x) = 
\lim_{\gamma_0\to 0} G_{\gamma_0, \mu_0, \sigma_0}(x) = \exp \Big (-\exp \Big (-\frac{x-\mu_0}{\sigma_0}\Big )\Big ).
\]
Further, we recall the concept of alpha (or strong) mixing as given in \cite{Bradley2005}. For this purpose, define the $\sigma$-fields
\[
    \mathcal{F}_{-\infty}^k = \sigma( \eps_i : i \leq k)~,~~
    \mathcal{F}_{j}^{\infty} = \sigma( \eps_i : i \geq j).
\]
and the strong mixing coefficients by
\begin{equation*}
    \alpha(n) = \sup_{A \in \mathcal{F}_{-\infty}^k, B \in \mathcal{F}_{j}^{\infty}, k < j - n} \left| P(A \cap B) - P(A)P(B) \right|.
\end{equation*}
The process $\{\eps_i\}_{i\in \mathbb{N} }$ is said to be \textbf{$\alpha$-mixing} (or strongly mixing) if:
\begin{equation*}
    \lim_{n \to \infty} \alpha(n) = 0.
\end{equation*}
With these preparations, we make the following assumptions.

\begin{assumption}\label{assump:kern}
	The kernel $K:\R\to\R$ is non-negative, supported on $[-1, 0]$ and satisfies $\int_{-1}^0K(x)\diff x = 1$. Further, it is Lipschitz-continuous on $(-1, 0)$.
\end{assumption}

\begin{assumption}\label{assump:mu}
	The function $\mu$ is twice differentiable with Lipschitz continuous second derivative. 
\end{assumption}

\begin{assumption}\label{assump:errors}
	For the (strictly) stationary error process $(\eps_i)_{i\in\Z}$, the following conditions hold:
	\begin{enumerate}
		\item[(i)] The error process $(\eps_i)_{i\in\Z}$ is centered, strongly mixing and satisfies
        \begin{align*}
            \mathbb{E}[|\eps_i|^{4+\delta}]&\leq C\\
            \sum_{i=1}^\infty \alpha(i)^{1/6-1/\delta}&<\infty
        \end{align*}
        for some $\delta>0$. We denote the long-run variance of $(\eps_i)_{i\in\Z}$, by 
				\[ \sigma^2 = \sum_{i\in\Z} \cov\big(\eps_{1+i}, \eps_1 \Big)\]
				 and note that it exists.
		\item[(ii)] The marginal distribution of the process $(|\eps_i|)_{i\in\Z}$ is in the domain of attraction of some GEV distribution $\Gc$, i.\,e., for a sequence $(\tilde{\eps}_i)_{i\in\Z}$ of i.i.d. copies of $(\eps_i)_{i\in\Z}$ it holds 
		\[ \lim_{n\to\infty} \pr( \max_{i=1}^n  |\tilde{\eps}_i| \le a_n x + b_n )= \Gc(x) \]
		for some sequences $a_n > 0$ and $b_n\in\R$. Further, $\pr(\max_{i=1}^n |\eps_i| \le a_n x + b_n)$ converges for some $x$. Finally, the sequences satisfy $\frac{a_n}{a_r}=\big(\frac{n}{r}\big)^\gamma$ and $b_n=b_r+a_r\frac{(\tfrac{n}{r})^\gamma-1}{\gamma}$, for $r=o(n)$.		
	\end{enumerate}
\end{assumption}

\begin{remark}\label{rem:conv_max}
  {\rm   Assumptions \ref{assump:kern} and \ref{assump:mu} are  standard in the time series literature (see, e.\,g., \citealp{bucher2021, heinrichs2021})
	Assumption \ref{assump:errors} is relatively mild as it is satisfied by many distributions of practical interest. The assumption about the scaling sequences is required to scale estimators of the GEV distribution from blocks of size $r$ to blocks of size $n$. Assumption \ref{assump:errors} (i) implies a common extremal mixing-type condition, introduced by \cite{leadbetter1974}, referred to as $D(u_n)$ condition. Together with Assumption \ref{assump:errors} (ii), it yields convergence of the distribution of $\tfrac{\max_{i=1}^n  |\eps_i| - b_n}{a_n}$ to some GEV distribution, see, e.\,g., Theorem 10.4 in \cite{beirlant2006}. Assumption (i) may be replaced by other weak dependence concepts such as physical dependence or $L^p$-m approximability, in this case however it might be necessary to establish (or assume) that the $D(u_n)$ condition holds as it does not follow straightforwardly from the definition of these dependence concepts.
    }
\end{remark}

\subsection{Critical Values} \label{sec:point}

When sequentially monitoring a time series for outliers, we conduct multiple statistical tests. If the number of tests is finite, we might use some correction procedure, such as the Bonferoni correction or the Holm-Bonferroni method, to control the overall level of the tests. If we choose this approach when monitoring for an indefinite amount of time, the level of each test must quickly converge to 0, thus loosing power. Conversely, in many applications a certain rate of type I errors might be acceptable, and we decide to keep the level, and therefore the power, constant.

In the following, we propose a decision rule and according levels such that (i) the level per test is constant or (ii) the overall level is bounded. In case (i), we bound the probability of a type I error per $n$ tests, and have to tolerate type I errors in the long run, while in case (ii), the power decreases over time, potentially not detecting existing outliers. This trade-off is inevitable, yet gives the applicant a theoretically sound choice. The asymptotic properties of the proposed decision rule are mainly based on the following theorem. 

\begin{theorem} \label{thm:conv} ~~~
	\begin{enumerate}
		\item If Assumption \ref{assump:errors} is satisfied, then, there exists some constant $\rho\in[0, 1]$ such that
		\[\lim_{n\to\infty} \pr\Big(\max_{i=k+1}^{k+n} |\eps_i |\le a_n x + b_n\Big)=\Gc^\rho(x)\]
		for any $k\in\N$, where $\Gc$ denotes the distribution from Assumption \ref{assump:errors}\,(ii).
		\item If Assumptions \ref{assump:kern}, \ref{assump:mu} and \ref{assump:errors} are satisfied, such that
		\begin{equation}\label{eq:conv_rate}
			\lim_{n\to\infty} \frac{1}{a_n}\bigg( h_n^3+ \sqrt{\frac{|\log(h_n)|}{nh_n}}+\frac{1}{nh_n^{5/4}}\bigg) = 0,
		\end{equation}
		and if $c_i = 0$ for $i = k+1, \dots, k+n$, for some $k\in\N$, then
		\[\lim_{n\to\infty} \pr\Big(\max_{i=k+1}^{k+n} |X_i - \tilde{\mu}_{n}(\tfrac{i}{n})|\le a_n x + b_n\Big)=\Gc^\rho(x).\]
	\end{enumerate}
\end{theorem}
We note that the limit $\Gc^\rho$ in the above theorem is a GEV distribution and denote its parameters by $\theta_\rho=(\gamma_\rho, \mu_\rho, \sigma_\rho)$. As the sequences $a_n$ and $b_n$ are unknown as well, we define $\theta_n = (\gamma_\rho, a_n \mu_\rho + b_n, a_n \sigma_\rho)$. Different approaches exist to approximate the unknown values of $\theta_n$, such as maximum likelihood estimation or the probability weighted moment (PWM) method. Denote an appropriate estimator of $\theta_n$ by $\hat{\theta}_n$.

\begin{assumption} \label{assump:estimation}
	The first observations $X_1, \dots, X_n$ do not contain outliers, i.\,e., $c_i=0$ for $i=1,\dots, n$. Let $\hat{\theta}_n=(\hat{\gamma}_n, \hat{\mu}_n, \hat{\sigma}_n)$ denote an estimator of $\theta_n$ based on these first $n$ observations such that $\hat{\gamma}_n - \gamma_\rho = o_\pr(1), \hat{\mu}_n - (a_n\mu_\rho + b_n) = o_\pr(a_n)$ and $\hat{\sigma}-a_n\sigma_\rho = o_\pr(a_n)$.
\end{assumption}

\newpage

\begin{remark} \label{rem:estimation} ~~~
{\rm 
\begin{enumerate}
    \item Let $q_{1-\alpha}(\theta)$ denote the $1-\alpha$ quantile of the GEV distribution with parameters $\theta$. Under Assumption \ref{assump:estimation}, straightforward calculations yield $\frac{q_{1-\alpha}(\hat{\theta}_n) - b_n}{a_n} = q_{1-\alpha}(\theta_\rho)+o_\pr(1)$.
	
	Assumption \ref{assump:estimation} is rather mild. For estimators based on maximum likelihood estimation or the probability weighted moment method, stronger results like asymptotic normality or almost sure convergence were shown under various conditions \citep[see, e.\,g.,][among others]{dombry2015, bucher2017, dombry2019, bucher2023}.

    \item Based on Assumption \ref{assump:estimation}, the role of the ``resolution'' $n$ with which $X_i$ is measured, is twofold. First, it is the length of the initial period without outliers, that is used to estimate the GEV distribution. Second, the proposed testing procedure, based on \eqref{eq:test}, bounds the joint level of $n$ consecutive tests. 
    
    Generally, the length of the initial period and the number of sequential tests with joint level $\alpha$, do not need to coincide, as long as they jointly grow to infinity. For the sake of readability, we assume that they match.
\end{enumerate}
}

\end{remark}
As argued before, for online outlier detection without finite time horizon, there is a trade-off between a loss of power, if the overall probability of a type I error $\alpha$ is fixed, and an accumulation of type I errors, if the power of each test should not vanish. Generally, we propose to reject $H_0^{(i)}$, if
\begin{equation} \label{eq:test}
	|X_i - \tilde{\mu}_{n}(\tfrac{i}{n})| > q_{1-\alpha_i}(\hat{\theta}_n).
\end{equation}
where $q_{1-\alpha}(\hat{\theta}_n)$ denotes the $1-\alpha$ quantile of the GEV distribution with parameters $\hat{\theta}_n$.

Depending on the application, we might choose to have a constant level $\alpha_i = \alpha$ for all time points $i\in\N$, possibly accumulating type I errors, or to select the level $\alpha_i = \bar{\alpha}_k$ for time points $i=(k-1) n+1,\dots, k n$, with $\sum_{k=1}^\infty \bar{\alpha}_k = \alpha$.  In the first case, the power remains constant over time, whereas in the second case it vanishes for large values of $i$.

The following corollary contains statements about the asymptotic level and consistency of the decision rule in \eqref{eq:test}. Therefore, based on $H_0^{(i)}: X_i \sim F_{\eps_i}\big(\cdot - \mu(\tfrac{i}{n})\big)$, define the joint null hypotheses 
\[
    H_0^{(k, \ell)} = \cap_{i=k}^{\ell} H_0^{(i)} ~\text{ and  } ~ \bar{H}_0 = \cap_{i=1}^\infty H_0^{(i)}.
\]
Note that under $\bar{H}_0$, it holds $X_i=\mu(\tfrac{i}{n})+ \eps_i$, for $i\in\N$. Even in the simplest case of independent, standard normally distributed errors, the quantity $\max_{i=1}^n |\eps_i|$ grows to infinity, as $n\to\infty$. To reliably detect outliers, the term $|c_i|$ must be sufficiently large, which is specified in the statement about the test's consistency.

\begin{corollary} \label{thm:test_properties}
	Let Assumptions \ref{assump:kern}, \ref{assump:mu}, \ref{assump:errors} and \ref{assump:estimation} be satisfied such that \eqref{eq:conv_rate} holds and let $\alpha > 0$.
	\begin{enumerate}
		\item[(i)] With constant level $\alpha_i = \alpha$, the decision rule in \eqref{eq:test} has asymptotic level $\alpha$ for $n$ consecutive tests, i.\,e., for any $k\ge n$,  \[\lim_{n\to\infty}\pr(\text{Reject}~H_0^{(i)}~\text{for some}~i\in \{k+1,\dots,k+n\}|H_0^{(k+1, k+n)})= \alpha.\]
		\item[(ii)] With $\alpha_i = \bar{\alpha}_{\lfloor i/n \rfloor}$, for $i\in\N$ and $\sum_{k=1}^\infty \bar{\alpha}_k = \alpha$, the decision rule in \eqref{eq:test} has asymptotic level $\alpha$, i.\,e., 
		\[\lim_{n\to\infty}\pr(\text{Reject}~H_0^{(i)}~\text{for some}~i\ge n|\bar{H}_0)= \alpha.\]
		\item[(iii)] The decision rule in \eqref{eq:test} is consistent against alternatives with $c_i=c_{i, n}$, such that $\pr(c_{i, n}\ge d_n)=1$, for $d_n \to\infty$ and 
		$\frac{\max\{|a_n|, |b_n|\}}{d_n}\to 0$, as $n\to\infty$. More specifically, 
		\[\lim_{n\to\infty}\pr(\text{Reject}~H_0^{(i)} | H_1^{(i)})= 1.\]
	\end{enumerate}
\end{corollary}

%% file: sec_empirical_results.tex
\section{Empirical Results} \label{sec:empirical}

In the following, we investigate the finite sample properties of the proposed sequential testing procedure and compare it with several  alternative procedures, which have been recently proposed in the literature. More precisely, we first use the decision rule proposed in \eqref{eq:test}, where the estimator $\tilde{\mu}_n(\tfrac{i}{n})$ is based on all observations until time $i$ and refer to it as the ``full'' version. Further, we use an adapted version of $\tilde{\mu}_n(\tfrac{i}{n})$, which is only calculated from the observations until time $i$, that are not detected as outliers. We will refer to this adapted version of the decision rule \eqref{eq:test} as the ``partial'' version. The rationale behind the partial version is that outliers potentially distort the estimator $\tilde{\mu}_n$, yielding larger deviations of $|X_i-\tilde{\mu}_n(\tfrac{i}{n})|$, so that the null hypothesis tends to be rejected more often, which impacts the test's finite sample properties.

Regarding the alternative methods, we use the testing procedures proposed by \cite{campulova2018} and \cite{holevsovsky2018} under the location-scale model. Further, we consider  an online version of the ARIMA model and reject the null hypothesis for large deviations of the residuals, as proposed by \cite{wette2024}. Finally, we use two established methods based on deep learning for the comparison  \citep{malhotra2015, munir2018}. In the following, we will refer to the respective methods as Ca2018, Ho2018, We2024, Ma2015 and Mu2018.
The first two methods are implemented in the \texttt{R} package \texttt{envoutliers} \citep{campulova2022}. Note that they are retrospective and cannot be used in an ``online'' fashion. The third method is implemented in the \texttt{Python} package \texttt{riverr} \citep{montiel2021}. For the latter two methods, we used custom \texttt{TensorFlow} implementations based on the model descriptions. The code for the subsequent experiments can be found in the GitHub repository: \url{https://github.com/FlorianHeinrichs/outlier_detection}.

For those methods, that do not control the joint error rates of multiple subsequent tests, we use the Bonferroni correction, which yields the individual level $\tfrac{\alpha}{n}$. Further, for all of these models, we used Chebychev's inequality for the calculation of critical values.

In Section \ref{sec:simulation_study} we conduct an extensive comparison by means of a simulation study and in Section \ref{sec:real_data}, we compare the methods when applied to real data.

\subsection{Simulation Study} \label{sec:simulation_study}

\begin{figure}
    \centering
    \includegraphics[width=0.9\linewidth]{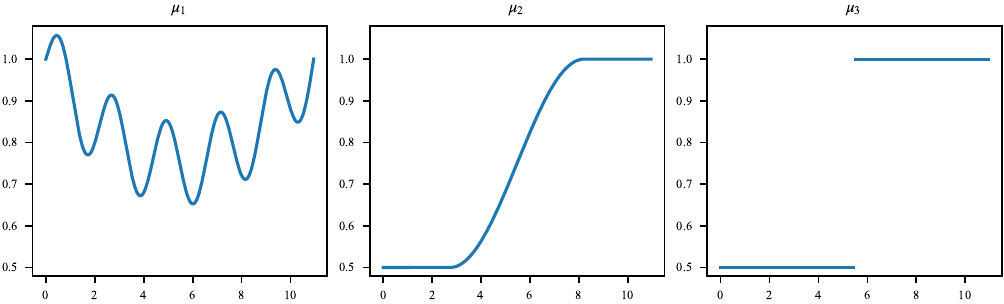}
    \caption{Different choices of the mean function $\mu$. Left: $\mu_1$. Center: $\mu_2$, Right: $\mu_3$.}
    \label{fig:mu}
\end{figure}

Given the model $X_i = \mu\big(\tfrac{i}{n}\big) + \eps_{i} + c_i$, as introduced in \eqref{eq:model}, the methods are compared with different choices of the mean function $\mu$, errors $\eps$ and levels of contamination $c$. For $\mu$, the choices
\begin{align*}
    \mu_0(t) & = 1, \quad & \mu_1(t) & = \big(\tfrac{t}{11}-\tfrac{1}{2}\big)^ 2 + \tfrac{1}{10}\sin(2 \pi  \tfrac{t}{11}) + \tfrac{3}{4}, \\
    \mu_2(t) & = \left\{ \begin{array}{ll}
        \tfrac{1}{2} &\text{if}~ t \le \tfrac{11}{4} \\
        \tfrac{3}{4} - \tfrac{1}{4}\sin(2\pi \tfrac{t}{11}) &\text{if}~ \tfrac{11}{4} < t < \tfrac{33}{4} \\
        1 &\text{if}~ t \ge \tfrac{33}{4},
    \end{array} \right.
    \quad & \mu_3(t) & = \left\{ \begin{array}{ll}
        \tfrac{1}{2} & \text{if}~t \le \tfrac{11}{2} \\
        1 &\text{if}~ t > \tfrac{11}{2},
    \end{array} \right.\\
\end{align*}
for $t\in [0, \infty)$, are studied, as displayed in Figure \ref{fig:mu}. The different mean functions were selected to cover continuous and non-continuous, monotone and non-monotone functions.

For the error process $\eps$, different distributions were considered. Recall that the cumulative distribution function of the Pareto distribution with parameter $\alpha$ is given by $1-(1+x)^{-\alpha}$. We considered independent random variables $(\eta_i)_{i\in\N}$ with
\begin{enumerate}
    \item[($\Nc$)  ] $\eta_i \sim \Nc(0, 1)$ (standard normal distribution)
    \item[($\Uc$)  ] $\eta_i \sim \Uc_{[0, 1]}$ (uniform distribution on $[0 ,1]$)
    \item[($Exp$)  ] $\eta_i \sim Exp(1)$ (exponential distribution with parameter $1$)
    \item[($Par_1$)] $\eta_i \sim Pareto(4)$ (Pareto distribution with parameter $4$)
    \item[($Par_2$)] $\eta_i \sim Pareto(2)$ (Pareto distribution with parameter $2$).
\end{enumerate}

Note that the variance in the last case is $\infty$. The random variables $\eta_i$ were used to construct (dependent) error processes

\[
\text{(IID)}~ \eps_i = \eta_i, \quad 
\text{(MA)}~ \eps_i = \eta_i + \tfrac{1}{2} \eta_{i-1}, \quad 
\text{(AR)}~ \eps_i = \eta_i + \tfrac{1}{2} \eps_{i-1}
\]
for $i\in \N$. Finally, we rescaled $\eps$ to have mean $0$ and standard deviation $\tfrac{1}{20}$, in all cases where the variance exists. The error processes were selected to include different dependence structures, and distributions in the domain of attraction of Gumbel, Fréchet and Weibull distributions.

We conducted experiments for $n=50, 100$ and $200$ and time horizons of length $11n$, i.\,e., $550, 1100$ and $2200$. Further, we considered time series with and without outliers. In the former case, outliers were randomly added to the original time series at 5\% of the observations. For each distribution, a minimal outlier height was defined in accordance with the assumption from Corollary \ref{thm:test_properties} (iii). This minimal height was multiplied with a random factor, sampled uniformly from $[-2, -1]\cup [1, 2]$. For each of the 180 combinations of mean function and error process, we generated $1\,000$ time series with outliers and another $1\,000$ time series without outliers.

For our method, we selected the one-sided quartic kernel $K(x) = \tfrac{15}{8}(1-x^2)^2$ with support on $[-1,0]$ and tuned the bandwidth with $5$-fold cross validation, minimizing the mean squared error. Finally, the joint level $\alpha$ for $n$ consecutive tests was set to $1\%$. 

In Tables \ref{table:iid_no_outliers} and \ref{table:iid_outliers}, we report the empirical specificity and sensitivity of the compared tests, for the experiments with i.i.d. data. Tables with results from all experiments can be found in Appendix \ref{app:appendix_empirical_results}. There were only minor differences between (IID), (MA) and (AR) errors, which are discussed below. Further, the deep learning based methods yielded basically trivial results as time series of length $n = 50, 100$ and $200$ are not sufficient to train neural networks reliably. Note that in the setting, that we study, only a single time series is given, thus one model per time series was trained. For the sake of readability, the results corresponding to the methods by \cite{malhotra2015} and \cite{munir2018} are only reported in the appendix.

As can be seen from Table \ref{table:iid_no_outliers}, in case of no outliers, the specificity is similarly high across all compared methods and increases with a growing sample size $n$. A notable exception is Ho2018, which has a poor specificity across all settings and even deteriorates as more data becomes available. The specificity is consistent across the different error distributions and mean functions. The two methods proposed in this paper, based on the full and partial time series, yield similar results, with the ``full'' version performing slightly better.

In the case of outliers, as reported in Table \ref{table:iid_outliers}, the results are generally consistent across different mean functions and tend to be better for larger sample sizes. While the specificity of Ho2018 is higher compared to the case of no outliers, it is still not competitive when compared to the remaining methods. As expected, the specificity of our ``full'' method is lower than that of the ``partial'' method. The observed specificities vary across the different error distributions, which can be explained when taking the respective sensitivities into account. More precisely, for uniformly and Pareto distributed errors, the detection of outliers seems to be more difficult and the null hypothesis seems to be accepted more often, compared to the other distributions. 

The sensitivity of We2024 is almost trivial across all distributions and mean functions. Limiting the joint level of $n$ successive tests seems to make the testing procedure overly conservative. Again, the sensitivity is similar for different mean functions and tends to be higher as the sample size $n$ increases. The differences between our ``full'' and ``partial'' methods seems negligible, which contrasts the aforementioned differences in specificity. The sensitivity of Ho2018 is almost $100\%$ across all cases, suggesting that the low specificity is a result of rejecting the null hypothesis frequently. For $n=50$, the sensitivity of Ca2018 is similar to that of the proposed methods, when the errors follow a normal or exponential distribution, and is larger for the other distributions. For $n=100$ and $n=200$, the proposed methods outperform Ca2018 for the normal, uniform and exponential distribution, while Ca2018 remains superior for the Pareto distribution. Generally, Ca2018 seems unable to detect outliers reliably, when the errors are uniformly distributed. Note that this is the only setting, where the data is in the domain of attraction of the Weibull distribution.

Regarding other error processes, the results generally coincide, with some notable exceptions. First, the specificity of our ``partial'' method is higher for independent errors compared to AR errors and slightly higher compared to MA errors. For the latter error processes, an extreme non-outlier value can propagate through time and lead to multiple rejected null hypotheses. For small sample sizes, the ``partial'' method's sensitivity is higher for uniform errors, the stronger their dependence. Contrarily, the sensitivity of the ``full'' method decreases with a stronger dependence.

In summary, the ML-based methods We2024, Ma2015 and Mu2018 seem to be overly conservative with sensitivities substantially smaller than $5\%$. This issue likely results from the Bonferroni correction, which is used to bound the familywise error rate of $n$ sequential tests to $\alpha$. Further, the neural networks cannot be trained well in this setting with low samples sizes, ranging from $50$ to $200$. A possible remedy for the difficulty of training neural networks with short time series, might be the use of functional neural networks that have fewer parameters and are known to work reasonably well, when data is sparse \citep{heinrichs2023}.
Generally, the proposed testing procedures and Ca2018 yield similarly accurate results, where the higher sensitivity of the former balances out the higher specificity of the latter. Note, however, that Ca2018 is a fully retrospective method and cannot be used in an online fashion. Further, our ``partial'' method seems to be slightly favorable compared to the ``full'' version. 

\setlength{\tabcolsep}{2pt}

\begin{table} 
\caption{Empirical specificity of various methods for IID errors, under different mean functions and error distributions.}
\begin{footnotesize}
\begin{tabular}{l| rrrrr | rrrrr | rrrrr}
	\hline \hline
	$n$ & \multicolumn{5}{c|}{50} & \multicolumn{5}{c|}{100} & \multicolumn{5}{c}{200} \\
	 $\mu$& $\Nc$ & $\Uc$ & $Exp$ & $Par_1$ & $Par_2$& $\Nc$ & $\Uc$ & $Exp$ & $Par_1$ & $Par_2$& $\Nc$ & $\Uc$ & $Exp$ & $Par_1$ & $Par_2$\\
	\hline
	\addlinespace[.2cm]
	\multicolumn{16}{l}{\quad\textit{Panel A: Specificity for $\mu_1$ }} \\
	Ours (full)& 98.3 & 98.5 & 98.2 & 97.6 & 97.5 & 99.3 & 99.3 & 99.4 & 99.2 & 99.1 & 99.7 & 99.7 & 99.8 & 99.8 & 99.8 \\
	Ours (partial)& 97.9 & 98.3 & 96.8 & 95.4 & 95.8 & 99.2 & 99.2 & 99.5 & 99.2 & 98.7 & 99.6 & 99.7 & 99.9 & 99.9 & 99.8 \\
	Ca2018& 99.4 & 100.0 & 99.4 & 99.0 & 98.6 & 99.6 & 100.0 & 99.7 & 99.3 & 98.8 & 99.8 & 100.0 & 99.8 & 99.4 & 98.9 \\
	Ho2018& 57.8 & 59.0 & 45.7 & 33.2 & 17.5 & 0.6 & 10.8 & 0.7 & 1.2 & 1.1 & 18.7 & 33.6 & 9.4 & 6.5 & 3.8 \\
	We2024& 99.7 & 99.7 & 99.7 & 99.7 & 99.6 & 99.8 & 99.8 & 99.8 & 99.8 & 99.8 & 99.9 & 99.9 & 99.9 & 99.9 & 99.9 \\
	\addlinespace[.2cm]
	\multicolumn{16}{l}{\quad\textit{Panel B: Specificity for $\mu_2$ }} \\
	Ours (full)& 98.6 & 98.6 & 98.6 & 98.2 & 98.1 & 99.3 & 99.3 & 99.6 & 99.5 & 99.3 & 99.7 & 99.7 & 99.8 & 99.8 & 99.7 \\
	Ours (partial)& 98.3 & 98.3 & 98.5 & 97.5 & 97.3 & 99.3 & 99.2 & 99.6 & 99.5 & 99.1 & 99.7 & 99.6 & 99.9 & 99.9 & 99.8 \\
	Ca2018& 99.4 & 100.0 & 99.4 & 99.0 & 98.6 & 99.6 & 100.0 & 99.6 & 99.3 & 98.8 & 99.8 & 100.0 & 99.8 & 99.4 & 99.0 \\
	Ho2018& 57.0 & 60.6 & 46.6 & 34.2 & 16.5 & 0.6 & 5.3 & 1.0 & 1.7 & 1.0 & 18.3 & 29.7 & 11.9 & 8.7 & 3.5 \\
	We2024& 99.7 & 99.7 & 99.7 & 99.7 & 99.6 & 99.8 & 99.8 & 99.8 & 99.8 & 99.8 & 99.9 & 99.9 & 99.9 & 99.9 & 99.9 \\
	\addlinespace[.2cm]
	\multicolumn{16}{l}{\quad\textit{Panel C: Specificity for $\mu_3$ }} \\
	Ours (full)& 97.5 & 97.6 & 98.1 & 97.8 & 97.8 & 98.9 & 98.8 & 99.4 & 99.3 & 99.3 & 99.5 & 99.5 & 99.7 & 99.8 & 99.7 \\
	Ours (partial)& 97.0 & 96.7 & 97.9 & 97.4 & 97.2 & 98.7 & 98.6 & 99.3 & 99.3 & 99.0 & 99.5 & 99.4 & 99.7 & 99.8 & 99.8 \\
	Ca2018& 99.6 & 99.6 & 98.7 & 98.5 & 98.4 & 99.7 & 99.7 & 98.9 & 98.9 & 98.7 & 99.7 & 99.8 & 99.3 & 99.3 & 98.9 \\
	Ho2018& 62.7 & 68.7 & 56.7 & 49.1 & 29.1 & 10.4 & 52.9 & 3.4 & 3.3 & 2.6 & 32.8 & 59.6 & 19.9 & 15.8 & 10.2 \\
	We2024& 99.7 & 99.7 & 99.7 & 99.7 & 99.6 & 99.8 & 99.8 & 99.8 & 99.8 & 99.8 & 99.9 & 99.9 & 99.9 & 99.9 & 99.9 \\
	\addlinespace[.2cm]
	\multicolumn{16}{l}{\quad\textit{Panel D: Specificity for $\mu_0$ }} \\
	Ours (full)& 98.3 & 98.6 & 98.8 & 98.5 & 97.9 & 99.3 & 99.3 & 99.6 & 99.5 & 99.4 & 99.7 & 99.7 & 99.8 & 99.8 & 99.8 \\
	Ours (partial)& 98.1 & 98.3 & 98.3 & 97.5 & 97.3 & 99.2 & 99.2 & 99.6 & 99.6 & 99.2 & 99.7 & 99.7 & 99.8 & 99.9 & 99.8 \\
	Ca2018& 99.5 & 100.0 & 99.3 & 99.0 & 98.6 & 99.7 & 100.0 & 99.7 & 99.3 & 98.8 & 99.8 & 100.0 & 99.8 & 99.4 & 99.0 \\
	Ho2018& 57.0 & 59.4 & 50.6 & 40.6 & 21.6 & 0.5 & 1.8 & 1.2 & 2.3 & 1.5 & 17.9 & 21.1 & 11.6 & 11.4 & 5.1 \\
	We2024& 99.7 & 99.7 & 99.7 & 99.6 & 99.6 & 99.8 & 99.8 & 99.8 & 99.8 & 99.8 & 99.9 & 99.9 & 99.9 & 99.9 & 99.9 \\
	\hline \hline 
 \end{tabular}
\end{footnotesize}
\label{table:iid_no_outliers}
\end{table}

\renewcommand{\arraystretch}{0.9}
\begin{table} 
\caption{Empirical specificity and sensitivity of various methods for IID errors, under different mean functions and error distributions.}
\begin{footnotesize}
\begin{tabular}{l| rrrrr | rrrrr | rrrrr}
	\hline \hline
	$n$ & \multicolumn{5}{c|}{50} & \multicolumn{5}{c|}{100} & \multicolumn{5}{c}{200} \\
	 $\mu$& $\Nc$ & $\Uc$ & $Exp$ & $Par_1$ & $Par_2$& $\Nc$ & $\Uc$ & $Exp$ & $Par_1$ & $Par_2$& $\Nc$ & $\Uc$ & $Exp$ & $Par_1$ & $Par_2$\\
	\hline
	\addlinespace[.2cm]
	\multicolumn{16}{l}{\quad\textit{Panel A: Specificity for $\mu_1$ }} \\
	Ours (full)& 86.8 & 95.2 & 85.5 & 89.4 & 88.3 & 87.4 & 96.7 & 87.6 & 93.4 & 92.6 & 87.2 & 97.3 & 88.0 & 95.7 & 94.2 \\
	Ours (partial)& 90.2 & 94.9 & 89.5 & 87.3 & 91.3 & 98.3 & 97.9 & 98.3 & 98.0 & 98.7 & 99.6 & 98.8 & 99.7 & 99.6 & 99.1 \\
	Ca2018& 99.9 & 100.0 & 99.9 & 99.7 & 99.7 & 100.0 & 100.0 & 99.9 & 99.8 & 99.8 & 100.0 & 100.0 & 99.9 & 99.8 & 99.9 \\
	Ho2018& 72.5 & 72.2 & 70.5 & 57.2 & 60.8 & 67.9 & 66.6 & 66.5 & 43.6 & 62.3 & 71.3 & 71.5 & 71.2 & 61.8 & 69.1 \\
	We2024& 99.7 & 99.7 & 99.7 & 99.7 & 99.7 & 99.9 & 99.8 & 99.8 & 99.8 & 99.8 & 99.9 & 99.9 & 99.9 & 99.9 & 99.9 \\
	\addlinespace[.2cm]
	\multicolumn{16}{l}{\quad\textit{Panel B: Specificity for $\mu_2$ }} \\
	Ours (full)& 87.2 & 95.7 & 87.8 & 92.5 & 90.2 & 87.5 & 96.8 & 88.5 & 94.9 & 92.9 & 87.3 & 97.4 & 88.1 & 95.9 & 94.3 \\
	Ours (partial)& 95.5 & 96.3 & 95.6 & 94.9 & 94.7 & 99.0 & 98.2 & 99.2 & 98.9 & 98.3 & 99.7 & 98.9 & 99.5 & 99.5 & 99.2 \\
	Ca2018& 99.9 & 100.0 & 99.9 & 99.7 & 99.7 & 100.0 & 100.0 & 99.9 & 99.8 & 99.8 & 100.0 & 100.0 & 99.9 & 99.9 & 99.9 \\
	Ho2018& 73.0 & 74.8 & 71.6 & 56.7 & 57.5 & 68.8 & 70.8 & 67.9 & 41.0 & 58.9 & 71.8 & 73.7 & 72.3 & 62.0 & 68.2 \\
	We2024& 99.7 & 99.7 & 99.7 & 99.7 & 99.7 & 99.8 & 99.8 & 99.8 & 99.8 & 99.8 & 99.9 & 99.9 & 99.9 & 99.9 & 99.9 \\
	\addlinespace[.2cm]
	\multicolumn{16}{l}{\quad\textit{Panel C: Specificity for $\mu_3$ }} \\
	Ours (full)& 86.8 & 94.7 & 87.3 & 91.9 & 90.1 & 87.1 & 96.4 & 88.4 & 94.8 & 92.8 & 87.0 & 97.1 & 87.8 & 96.0 & 93.8 \\
	Ours (partial)& 95.0 & 95.2 & 95.0 & 95.3 & 95.4 & 98.6 & 97.7 & 98.3 & 98.9 & 98.3 & 99.4 & 98.5 & 99.4 & 99.5 & 99.3 \\
	Ca2018& 99.9 & 99.7 & 99.8 & 99.7 & 99.7 & 99.9 & 99.9 & 99.9 & 99.8 & 99.8 & 100.0 & 99.9 & 99.9 & 99.8 & 99.9 \\
	Ho2018& 67.8 & 70.6 & 63.6 & 42.5 & 53.3 & 63.7 & 65.8 & 60.8 & 23.5 & 57.1 & 69.1 & 70.8 & 68.0 & 51.5 & 65.0 \\
	We2024& 99.7 & 99.7 & 99.7 & 99.7 & 99.7 & 99.9 & 99.9 & 99.8 & 99.8 & 99.8 & 99.9 & 99.9 & 99.9 & 99.9 & 99.9 \\
	\addlinespace[.2cm]
	\multicolumn{16}{l}{\quad\textit{Panel D: Specificity for $\mu_0$ }} \\
	Ours (full)& 87.6 & 96.0 & 87.4 & 92.3 & 90.6 & 87.7 & 96.8 & 88.5 & 94.9 & 92.9 & 87.4 & 97.3 & 88.1 & 95.5 & 93.7 \\
	Ours (partial)& 95.4 & 96.7 & 96.1 & 96.0 & 95.4 & 99.2 & 98.2 & 98.7 & 98.9 & 98.7 & 99.7 & 98.8 & 99.4 & 99.6 & 99.3 \\
	Ca2018& 99.9 & 100.0 & 99.9 & 99.7 & 99.7 & 100.0 & 100.0 & 99.9 & 99.8 & 99.8 & 100.0 & 100.0 & 99.9 & 99.9 & 99.9 \\
	Ho2018& 73.3 & 75.8 & 72.3 & 58.0 & 58.0 & 69.1 & 72.3 & 68.5 & 41.1 & 59.1 & 71.8 & 74.4 & 72.8 & 62.9 & 68.5 \\
	We2024& 99.7 & 99.7 & 99.7 & 99.7 & 99.7 & 99.8 & 99.9 & 99.9 & 99.8 & 99.8 & 99.9 & 99.9 & 99.9 & 99.9 & 99.9 \\
	\addlinespace[.2cm]
	\multicolumn{16}{l}{\quad\textit{Panel E: Sensitivity for $\mu_1$ }} \\
	Ours (full)& 80.4 & 32.7 & 79.1 & 52.4 & 63.8 & 93.8 & 48.7 & 91.3 & 60.9 & 61.2 & 97.5 & 61.1 & 95.4 & 60.6 & 58.0 \\
	Ours (partial)& 82.9 & 31.6 & 81.4 & 56.5 & 65.5 & 94.9 & 47.0 & 92.6 & 61.2 & 62.0 & 99.1 & 59.8 & 96.8 & 61.5 & 58.9 \\
	Ca2018& 82.9 & 41.5 & 82.9 & 78.3 & 79.4 & 86.1 & 38.5 & 85.7 & 82.6 & 84.0 & 87.2 & 40.3 & 87.1 & 84.0 & 87.2 \\
	Ho2018& 100.0 & 99.3 & 100.0 & 100.0 & 100.0 & 100.0 & 99.9 & 100.0 & 100.0 & 94.0 & 100.0 & 100.0 & 99.6 & 100.0 & 100.0 \\
	We2024& 2.0 & 0.5 & 3.1 & 2.4 & 3.1 & 1.1 & 0.2 & 1.6 & 1.3 & 1.7 & 0.6 & 0.1 & 0.9 & 0.7 & 0.9 \\
	\addlinespace[.2cm]
	\multicolumn{16}{l}{\quad\textit{Panel F: Sensitivity for $\mu_2$ }} \\
	Ours (full)& 90.6 & 49.3 & 87.0 & 64.6 & 63.9 & 96.4 & 56.1 & 92.6 & 59.6 & 61.1 & 97.6 & 60.6 & 95.2 & 60.1 & 58.2 \\
	Ours (partial)& 92.9 & 47.9 & 89.0 & 66.5 & 65.9 & 97.7 & 54.3 & 94.0 & 60.5 & 62.5 & 99.1 & 59.3 & 96.8 & 61.0 & 59.1 \\
	Ca2018& 84.1 & 45.4 & 84.0 & 79.6 & 80.9 & 86.4 & 46.7 & 86.3 & 82.6 & 85.0 & 87.4 & 46.5 & 87.4 & 84.7 & 87.8 \\
	Ho2018& 100.0 & 99.4 & 100.0 & 99.9 & 99.9 & 100.0 & 100.0 & 100.0 & 100.0 & 91.4 & 100.0 & 100.0 & 98.9 & 100.0 & 100.0 \\
	We2024& 2.1 & 0.5 & 3.1 & 2.4 & 3.0 & 1.1 & 0.2 & 1.7 & 1.3 & 1.7 & 0.6 & 0.1 & 0.8 & 0.7 & 0.9 \\
	\addlinespace[.2cm]
	\multicolumn{16}{l}{\quad\textit{Panel G: Sensitivity for $\mu_3$ }} \\
	Ours (full)& 90.7 & 49.0 & 86.5 & 62.9 & 64.5 & 96.5 & 56.3 & 91.7 & 60.9 & 61.7 & 97.5 & 60.6 & 95.5 & 58.4 & 60.4 \\
	Ours (partial)& 91.6 & 47.1 & 88.1 & 64.1 & 66.5 & 97.6 & 54.6 & 93.5 & 61.5 & 62.9 & 99.1 & 59.2 & 97.2 & 59.2 & 61.2 \\
	Ca2018& 77.3 & 33.2 & 81.3 & 74.1 & 80.0 & 83.5 & 30.6 & 85.1 & 80.8 & 84.3 & 86.0 & 27.9 & 86.9 & 83.6 & 87.3 \\
	Ho2018& 100.0 & 98.7 & 100.0 & 100.0 & 100.0 & 100.0 & 99.7 & 100.0 & 100.0 & 93.6 & 99.9 & 99.9 & 99.8 & 100.0 & 100.0 \\
	We2024& 2.1 & 0.5 & 3.1 & 2.4 & 3.1 & 1.1 & 0.2 & 1.6 & 1.3 & 1.7 & 0.6 & 0.1 & 0.9 & 0.7 & 0.9 \\
	\addlinespace[.2cm]
	\multicolumn{16}{l}{\quad\textit{Panel H: Sensitivity for $\mu_0$ }} \\
	Ours (full)& 90.8 & 48.7 & 87.7 & 63.2 & 63.9 & 96.1 & 55.7 & 91.5 & 61.0 & 59.5 & 97.5 & 61.5 & 94.6 & 61.2 & 59.7 \\
	Ours (partial)& 92.9 & 46.9 & 89.3 & 64.7 & 65.6 & 97.6 & 54.3 & 93.2 & 62.1 & 60.6 & 99.1 & 60.3 & 96.3 & 62.0 & 60.5 \\
	Ca2018& 84.6 & 46.4 & 84.4 & 80.1 & 80.7 & 86.9 & 50.6 & 86.5 & 82.6 & 85.1 & 87.5 & 52.6 & 87.7 & 84.8 & 87.7 \\
	Ho2018& 100.0 & 99.4 & 100.0 & 99.8 & 99.6 & 100.0 & 100.0 & 100.0 & 100.0 & 90.9 & 100.0 & 100.0 & 99.1 & 100.0 & 100.0 \\
	We2024& 2.1 & 0.4 & 3.1 & 2.3 & 3.2 & 1.1 & 0.2 & 1.7 & 1.3 & 1.7 & 0.6 & 0.1 & 0.9 & 0.7 & 0.9 \\
	\hline \hline 
 \end{tabular}
\end{footnotesize}
\label{table:iid_outliers}
\end{table}

\setlength{\tabcolsep}{6pt}
\renewcommand{\arraystretch}{1}

\subsection{Real Data Experiments} \label{sec:real_data}

Time series with non-stationary mean function naturally arise in the field of meteorology. We considered daily temperature data from eight Australian cities, spanning approximately 150 years \citep{australian_climate}\footnote{The data can be accessed freely via \url{http://www.bom.gov.au/climate/data/index.shtml}.}. This data can be assumed to contain a seasonal component and a trend. These time series have been studied in the literature, and the null hypothesis of a stationary mean can be rejected for all cities, but Gunnedah, at level $\alpha=5\%$, even when seasonality is taking into account \citep{bucher2020}. 

During the monitored periods of more than a century, measurement errors occurred. Sometimes temperatures were not measured (correctly) on individual days, and sometimes over longer periods of time. For periods of 3 or more consecutive days without a recorded temperature, we imputed the mean over all years of the corresponding day, i.\,e., 
\[ T_t = \frac{1}{|D_t|} \sum_{i \in D_t} X_i, \]
where $D_t$ denotes the set of observations of the same day and month as the time point $t$, but with non-missing values. These imputed values and correctly measured temperatures are considered as ``normal'' observations. 

\begin{figure}
    \centering
    \includegraphics[width=0.9\linewidth]{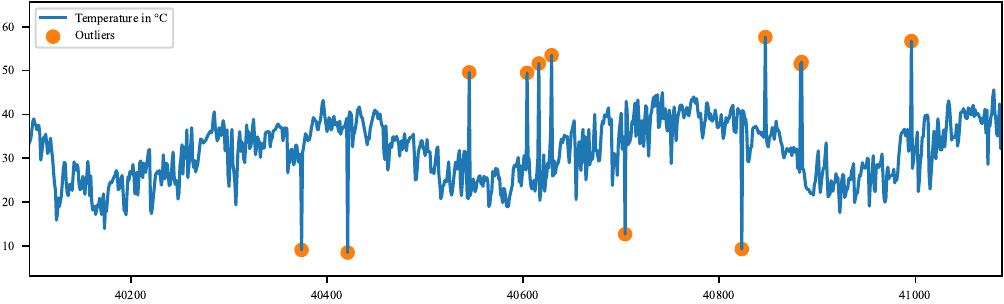} \\
    \includegraphics[width=0.9\linewidth]{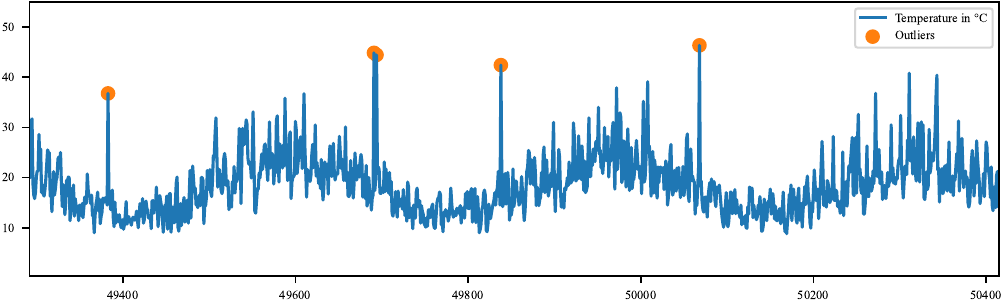} \\
    \caption{Temperature data from different Australian locations with marked outliers. Top: Boulia. Bottom: Hobart.}
    \label{fig:real_data}
\end{figure}

For single missing values and periods of two days, we imputed values that were considered as ``outliers''. More specifically, we calculated the ``moving range'' $r_t$ as the difference between maximum and minimum temperature over the last 365 days and the corresponding moving mean $m_t$. With $T_t$ as before, we imputed the value
\[ X_t = \left\{ \begin{array}{rr}
     T_t - r_t & \text{if}~ T_t > m_t \\
     T_t + r_t & \text{if}~ T_t \le m_t 
\end{array} \right. \]
for time instances $t$ with missing values. These imputed values simulate measurement errors in situations, where sensors record incorrect values instead of placeholders for ``missing values''. Excerpts of the processed data from Boulia (between days 40100 and 41100) and from Hobart (between days 49300 and 50400) of approximately 3 years are displayed in Figure \ref{fig:real_data}. 

For our experiments, we set $n=365$, which corresponds to an initial year without outliers, and chose $\alpha=1\%$, as before. Clearly, the first year does not need to be free of measurement errors. We selected the first consecutive 365 days without outlier as the initial period and discarded the data before. 

We used generally the same options for the different methods as in Section \ref{sec:simulation_study}. To select an adequate bandwidth for the Jackknife estimator of the mean function $\mu$, we used cross validation with bandwidths between 30 and 50 days. 

The results reflect those of the simulation study. The machine learning-based methods Ma2015, Mu2018 and We2024 do not detect any outliers in the data. Ho2018 detects all outliers, but has a (very) low specificity, which renders the method practically useless. As before, the specificities of the proposed testing procedures and Ca2018 are similar, while the sensitivities are (substantially) higher in all cities except Boulia and Sydney.

\setlength{\tabcolsep}{4pt}
\begin{table}
\caption{Empirical specificity and sensitivity of various methods for temperature data in different Australian locations.}
\begin{footnotesize}
\begin{tabular}{l| rrrrrrrr}
	\hline \hline
	City & Boulia & Gayndah & Gunnedah & Hobart & Melbourne & Cape Otway & Robe & Sydney\\
	\hline
	\addlinespace[.2cm]
	\multicolumn{9}{l}{\quad\textit{Panel A: Specificity}} \\
	Ours (full)& 99.9 & 99.5 & 99.0 & 97.1 & 100.0 & 99.6 & 99.4 & 100.0 \\
	Ours (partial)& 99.9 & 99.7 & 99.9 & 97.3 & 100.0 & 99.4 & 99.8 & 100.0 \\
	Ca2018& 99.9 & 100.0 & 99.9 & 99.9 & 99.1 & 99.7 & 99.5 & 99.9 \\
	Ho2018& 55.1 & 48.1 & 34.6 & 36.3 & 31.0 & 0.7 & 0.7 & 37.4 \\
	We2024& 100.0 & 100.0 & 100.0 & 100.0 & 100.0 & 100.0 & 100.0 & 100.0 \\
	Ma2015& 100.0 & 100.0 & 100.0 & 100.0 & 100.0 & 100.0 & 100.0 & 100.0 \\
	Mu2018& 100.0 & 100.0 & 100.0 & 100.0 & 100.0 & 100.0 & 100.0 & 100.0 \\
	\addlinespace[.2cm]
	\multicolumn{9}{l}{\quad\textit{Panel B: Sensitivity}} \\
	Ours (full)& 40.0 & 97.3 & 95.7 & 99.6 & 100.0 & 49.8 & 78.7 & 89.2 \\
	Ours (partial)& 45.3 & 100.0 & 98.4 & 100.0 & 100.0 & 60.8 & 96.0 & 92.6 \\
	Ca2018& 65.3 & 83.0 & 37.8 & 91.8 & 100.0 & 23.3 & 45.1 & 97.3 \\
	Ho2018& 100.0 & 100.0 & 100.0 & 100.0 & 100.0 & 100.0 & 100.0 & 100.0 \\
	We2024& 0.0 & 0.0 & 0.0 & 0.0 & 0.0 & 0.0 & 0.0 & 0.0 \\
	Ma2015& 0.0 & 0.0 & 0.0 & 0.0 & 0.0 & 0.0 & 0.0 & 0.0 \\
	Mu2018& 0.0 & 0.0 & 0.0 & 0.0 & 0.0 & 0.0 & 0.0 & 0.0 \\
	\addlinespace[.2cm]
	\multicolumn{9}{l}{\quad\textit{Panel C: Statistics}} \\
	No. data points & 49\,140 & 37\,809 & 36\,162 & 52\,163 & 58\,324 & 57\,323 & 50\,318 & 59\,047 \\
	No. outliers & 956 & 182 & 2\,256 & 233 & 2 & 2\,468 & 916 & 148 \\
	\hline \hline 
 \end{tabular}
\end{footnotesize}
\end{table}
\setlength{\tabcolsep}{6pt}

%% file: app_main_thm.tex
\subsection{Proof of Theorem \ref{thm:conv}} \label{app:main_thm}

The first part of the theorem follows directly from Theorem 10.4 in \cite{beirlant2006}. 
By assumption, $c_i=0$, so that $X_i = \mu(\tfrac{i}{n}) + \eps_i$, for $i\in\{k+1,\dots, k+n\}$. By the reverse triangle inequality, we have
\begin{align*}
	\bigg| \max_{i=k+1}^{k+n} |X_i - \tilde{\mu}_{n}(\tfrac{i}{n})| - \max_{i=k+1}^{k+n}|\eps_i|  \bigg|
	&\le \max_{i=k+1}^{k+n} |X_i - \tilde{\mu}_{n}(\tfrac{i}{n}) - \eps_i| = \max_{i=k+1}^{k+n} |\mu(\tfrac{i}{n}) - \tilde{\mu}_{n}(\tfrac{i}{n})|.
\end{align*}
In the following, we show that $\max_{i=k+1}^{k+n} |\mu(\tfrac{i}{n}) - \tilde{\mu}_{n}(\tfrac{i}{n})| = o_\pr(a_n)$. Together with the theorem's first part, this finalizes the proof.

By Proposition \ref{thm:lle}, it holds
\begin{equation} \label{eq:main_approx}
	\max_{i=k+1}^{k+n} |\mu(\tfrac{i}{n}) - \tilde{\mu}_{n}(\tfrac{i}{n})| 
	= \max_{i=k+1}^{k+n} \bigg| \frac{1}{nh_n}\sum_{j\in\N} \eps_j K^*\big(\tfrac{j-i}{nh_n}\big)  \bigg|
	+\Oc(h_n^3 + \tfrac{h_n}{n}) + \Oc_\pr\Big(\tfrac{1}{nh_n^{5/4}}\Big). 
\end{equation}
Recall that $\sigma^2$ denotes the long-run variance of the error process $(\eps_i)_{i\in\Z}$ and $\|\cdot\|_2$ the $L^2$-norm, i.\,e., $\|f\|_2^2=\int f^2(x)\diff x$. With this notation, define the process 
\[ Z_{n}(t) = \frac{1}{\sqrt{nh_n} \sigma \|K^*\|_2}\sum_{j\in\N} \eps_j K^*\big(\tfrac{j-nt}{nh_n}\big).\]
Note that $Z_n(t)$ is essentially the same process as $Z_{n, 1}(t)$ in the proof of Theorem A.2 in the appendix of \cite{bucher2021}. By analogous arguments, using Theorem from \cite{Shao1987} instead of proposition 5 from \cite{Zhou2013}, it follows that 
\[ \ell_n \sup_{t\in [k/n, k/n+1]}|Z_n(t)| - \ell_n^2 \convw G_{0, \log(2), 1}, \]
where $\ell_n = \sqrt{2 \log\big(\tfrac{\|(K^*)'\|_2}{\|K^*\|_2 2\pi h_n}\big)}$. Since $\ell_n=\Oc(\sqrt{|\log(h_n)|})$, it follows from \eqref{eq:main_approx} that
\[	\max_{i=k+1}^{k+n} |\mu(\tfrac{i}{n}) - \tilde{\mu}_{n}(\tfrac{i}{n})| = \Oc(h_n^3+\tfrac{h_n}{n}) + \Oc_\pr\Big(\sqrt{\tfrac{|\log(h_n)|}{nh_n}}+\tfrac{1}{nh_n^{5/4}}\Big),\]
where the right-hand side is of order $o_\pr(a_n)$ by \eqref{eq:conv_rate}.

%% file: app_test_properties.tex
\subsection{Proof of Corollary \ref{thm:test_properties}} \label{app:test_properties}

\begin{enumerate}
	\item[(i)] 
	Note that we can rewrite 
	\begin{align*}
		& \pr(\text{Reject}~H_0^{(i)}~\text{for some}~i\in \{k+1,\dots,k+n\}|H_0^{(k+1, k+n)}) \\
		& = \pr\Big(\max_{i=k+1}^{k+n} |X_i - \tilde{\mu}_n(\tfrac{i}{n})| > q_{1-\alpha}(\hat{\theta}_n) |H_0^{(k+1, k+n)}\Big).
	\end{align*}
	By Remark \ref{rem:estimation}, it holds $\frac{q_{1-\alpha}(\hat{\theta}_n) - b_n}{a_n} = q_{1-\alpha}(\theta_\rho)+o_\pr(1)$. Together with the second part of Theorem \ref{thm:conv}, (i) follows.
	
	\item[(ii)] Similarly, we rewrite
	\begin{align*}
		& \pr(\text{Reject}~H_0^{(i)}~\text{for some}~i\in \N|\bar{H}_0)\\
		& = \pr(\max_{i\in\N} |X_i - \tilde{\mu}_n(\tfrac{i}{n})| > q_{1-\alpha_i}(\hat{\theta}_n) |\bar{H}_0)\\
		& \le \sum_{k=1}^{\infty} \pr\Big(\max_{i=(k-1)n+1}^{kn} |X_i - \tilde{\mu}_n(\tfrac{i}{n})| > q_{1-\bar{\alpha}_{\lceil i/n \rceil}}(\hat{\theta}_n) \Big|\bar{H}_0\Big),
	\end{align*}
	which converges to $\sum_{k=1}^{\infty} \bar{\alpha}_k = \alpha$ as $n\to\infty$, by the same arguments as before.
	
	\item[(iii)] From Remark \ref{rem:estimation}, we obtain $q_{1-\alpha}(\hat{\theta}_n) = b_n + \Oc_\pr(a_n)$. By Proposition \ref{thm:lle}, it follows $|X_i - \tilde{\mu}_n(\tfrac{i}{n})| = |\mu(\tfrac{i}{n}) + \eps_i + c_i - \tilde{\mu}_n(\tfrac{i}{n})| = |c_i| + \Oc_\pr(1)$. Further note that 
	\[ \pr\big(\text{Reject}~H_0^{(i)} | H_1^{(i)}\big) = \pr\big( |X_i - \tilde{\mu}_n(\tfrac{i}{n})| > q_{1-\alpha}(\hat{\theta}_n) | H_1^{(i)}\big). \]	
	With $\frac{\max\{|a_n|, |b_n|\}}{d_n}\to 0$, the corollary's statement follows.
\end{enumerate}

%% file: app_jackknife_estimator.tex
\subsection{Jackknife Estimation} \label{app:appendix_jackknife}

For the Jackknife estimator $\tilde{\mu}_{n}$, define its corresponding kernel as 
\[K^*(x)= 2 \sqrt{2} \bar{K}(\sqrt{2}x) - \bar{K}(x),\] 
with $\bar{K}(x)=\frac{\kappa_2-\kappa_1 x}{\kappa_0\kappa_2-\kappa_1^2}K(x)$ and $\kappa_\ell = \int_{-1}^0 x^\ell K(x)\diff x$, for $\ell=0, 1, 2, 3$.

\begin{proposition}\label{thm:lle}
	Let $h_n\searrow 0$ as $n\to\infty$. If Assumptions \ref{assump:kern} and \ref{assump:mu} are satisfied, and the errors satisfy $\ex[\eps_i^4]<\infty$, then, for any $\delta > h_n$,
	\[\sup_{t \in [\delta, 1+\delta]} \bigg| \tilde{\mu}_{n}(t)- \mu(t) 
	- \frac{1}{nh_n}\sum_{i\in\N} \eps_i K^*\big(\tfrac{i-nt}{nh_n}\big)  \bigg|
	=\Oc(h_n^3 + \tfrac{h_n}{n}) + \Oc_\pr\Big(\tfrac{1}{nh_n^{5/4}}\Big).\]
	
\end{proposition}
To calculate $\hat{\mu}_{h_n}(t)$, we need to minimize
\[ g_{n}(b_0, b_1) = \frac{1}{nh_n}\sum_{i\in\N} \Big(X_{i} - b_0 - b_1 \big(\tfrac{i}{n}-t\big) \Big)^2 K\big(\tfrac{i-nt}{nh_n}\big). \]
By convexity of $g_{n}$, $(b_0, b_1)$ is a global minimum if the gradient of $g_{n}$ vanishes. Straightforward calculations yield
\[ \frac{\partial g_{n}}{\partial b_\ell} = - 2 \big(h_n^\ell R_{\ell}(t)-b_0h_n^\ell S_{\ell}(t)-b_1 h_n^{\ell+1}S_{\ell+1}(t)\big), \]
for $\ell=0, 1$, where
\[S_{\ell}(t) = \frac{1}{nh_n} \sum_{i\in\N} (\tfrac{i-nt}{nh_n}\big)^\ell K\big(\tfrac{i-nt}{nh_n}\big)
\quad\text{and}\quad
R_{\ell}(t) = \frac{1}{nh_n} \sum_{i\in\N} X_i (\tfrac{i-nt}{nh_n}\big)^\ell K\big(\tfrac{i-nt}{nh_n}\big).\]
The gradient of $g_{n}$ vanishes if and only if 
\begin{align*}
	\left(\begin{array}{c} b_0 \\ b_1 \end{array}\right)
	= \left(\begin{array}{cc}
		S_{0}(t) & h_n S_{1}(t)\\
		h_n S_{1}(t) & h_n^2 S_{2}(t)
	\end{array}\right)^{-1}  \left(\begin{array}{c}
		R_{0}(t)\\
		h_n R_{1}(t)\\
	\end{array}\right).
\end{align*}
By a Taylor expansion and Lipschitz continuity of $\mu''(t)$, it follows 
\begin{align}\label{eq:taylor_expansion}
	\sup_{t \ge h_n} \bigg| &\hat{\mu}_{n}(t)- \mu(t) - \tfrac{h_n^2}{2} \mu''(t) \tfrac{S_{2}^2(t)-S_{1}(t)S_{3}(t)}{S_{0}(t)S_{2}(t)-S_{1}^2(t)}\\
	&- \frac{1}{nh_n}\sum_{i\in\N} \tfrac{S_{2}(t)-S_{1}(t)\big(\frac{i-nt}{nh_n}\big)}{S_{0}(t)S_{2}(t)-S_{1}^2(t)}K\big(\tfrac{i-nt}{nh_n}\big)\eps_i  \bigg| =\Oc(h_n^3).\notag
\end{align}
By Lipschitz continuity of $K$, $S_{\ell}(t) = \kappa_\ell + \Oc(\tfrac{1}{nh_n})$, uniformly for $t\ge h_n$. In particular, it holds 
\begin{equation}\label{eq:conv_kern_const}
	\sup_{t \ge h_n}\bigg| \frac{S_{2}^2(t)-S_{1}(t)S_{3}(t)}{S_{0}(t)S_{2}(t)-S_{1}^2(t)} - \frac{\kappa_2^2-\kappa_1\kappa_3}{\kappa_0\kappa_2-\kappa_1^2} \bigg| = \Oc(\tfrac{1}{nh_n})
\end{equation}
and 
\begin{equation}\label{eq:conv_kern}
	\sup_{t \ge h_n}|\Delta_n(i, t)| = \sup_{t \ge h_n}\bigg| \frac{S_{2}(t)-S_{1}(t)\big(\frac{i-nt}{nh_n}\big)}{S_{0}(t)S_{2}(t)-S_{1}^2(t)}K\big(\tfrac{i-nt}{nh_n}\big) - \bar{K}\big(\tfrac{i-nt}{nh_n}\big)\bigg| = \Oc(\tfrac{1}{nh_n}), 
\end{equation}
where $\bar{K}(x)=\frac{\kappa_2-\kappa_1 x}{\kappa_0\kappa_2-\kappa_1^2}K(x)$ and $\Delta_n(i, t)$ is defined in an obvious manner. From \eqref{eq:taylor_expansion} and \eqref{eq:conv_kern_const}, we obtain the uniform approximation
\[ \sup_{t \ge h_n} \bigg| \hat{\mu}_{n}(t)- \mu(t) - \tfrac{h_n^2}{2} \mu''(t) \tfrac{\kappa_2^2-\kappa_1\kappa_3}{\kappa_0\kappa_2-\kappa_1^2}
- \frac{1}{nh_n}\sum_{i\in\N} \tfrac{S_{2}(t)-S_{1}(t)\big(\frac{i-nt}{nh_n}\big)}{S_{0}(t)S_{2}(t)-S_{1}^2(t)}K\big(\tfrac{i-nt}{nh_n}\big)\eps_i  \bigg| =\Oc\big(h_n^3 + \tfrac{h_n}{n}\big). \]
By Hölder's inequality, we obtain, for $\delta > h_n$,
\begin{equation}\label{eq:bound_moment} 
	\ex\bigg[ \sup_{t \in [\delta, 1+\delta]} \bigg( \sum_{i\in\N} \Delta_n(i, t)\eps_i \bigg)^4 \bigg] \le \sup_{t \in [\delta, 1+\delta]} \bigg( \sum_{i=\lceil (\delta-h_n)n\rceil}^{\lfloor (1+\delta) n\rfloor} \Delta_n^{4/3}(i, t)\bigg)^3 \sum_{i=\lceil (\delta-h_n)n\rceil}^{\lfloor (1+\delta) n\rfloor} \ex[ \eps_i^4].   
\end{equation}
By assumption, the second sum on the right-hand side of the previous display is of order $\Oc(n)$. Further, $\Delta_n(i, t)=0$ for $i \notin [n(t-h_n), nt]$. By \eqref{eq:conv_kern}, it exists some constant $C\in\R$ such that
\[ \sum_{i=\lceil (\delta-h_n)n\rceil}^{\lfloor (1+\delta) n\rfloor} \Delta_n^{4/3}(i, t) 
\le \frac{C}{(nh_n)^{1/3}} \]

uniformly in $t\in[\delta, 1+\delta]$. Thus, the right-hand side of \eqref{eq:bound_moment} is of order $\Oc(h_n^{-1})$ and
\[\sup_{t \in [\delta, 1+\delta]} \frac{1}{nh_n}\sum_{i\in\N} \Delta_n(i, t)\eps_i = \Oc_\pr\Big(\tfrac{1}{nh_n^{5/4}}\Big).\]
It follows from \eqref{eq:taylor_expansion} and \eqref{eq:conv_kern_const} that
\begin{align*}
	&\sup_{t \in [\delta, 1+\delta]} \bigg| \hat{\mu}_{h_n}(t)- \mu(t) - \tfrac{h_n^2}{2} \mu''(t) \frac{\kappa_2^2-\kappa_1\kappa_3}{\kappa_0\kappa_2-\kappa_1^2}
	- \frac{1}{nh_n}\sum_{i\in\N} \bar{K}\big(\tfrac{i-nt}{nh_n}\big)\eps_i  \bigg|\\
	& =\Oc(h_n^3+ \tfrac{h_n}{n}) + \Oc_\pr\Big(\tfrac{1}{nh_n^{5/4}}\Big).
\end{align*}
The statement of the theorem   follows from the definition of the Jackknife estimator.

\bigskip

\noindent
{\bf Acknowledgments}
This work was supported by TRR 391 \textit{Spatio-temporal Statistics for the Transition of Energy and Transport} (Project number 520388526) funded by the Deutsche Forschungsgemeinschaft (DFG, German Research Foundation).

%% file: app_empirical_results.tex
\section{Additional Empirical Results} \label{app:appendix_empirical_results}

\setlength{\tabcolsep}{3pt}
\renewcommand{\arraystretch}{0.85}

\begin{table}[h]
\caption{Empirical specificity of various methods for IID errors, under different mean functions and error distributions.}
\begin{footnotesize}
\begin{tabular}{l| rrrrr | rrrrr | rrrrr}
	\hline \hline
	$n$ & \multicolumn{5}{c|}{50} & \multicolumn{5}{c|}{100} & \multicolumn{5}{c}{200} \\
	 $\mu$& $\Nc$ & $\Uc$ & $Exp$ & $Par_1$ & $Par_2$& $\Nc$ & $\Uc$ & $Exp$ & $Par_1$ & $Par_2$& $\Nc$ & $\Uc$ & $Exp$ & $Par_1$ & $Par_2$\\
	\hline
	\addlinespace[.2cm]
	\multicolumn{16}{l}{\quad\textit{Panel A: Specificity for $\mu_1$ }} \\
	Ours (full)& 98.3 & 98.5 & 98.2 & 97.6 & 97.5 & 99.3 & 99.3 & 99.4 & 99.2 & 99.1 & 99.7 & 99.7 & 99.8 & 99.8 & 99.8 \\
	Ours (partial)& 97.9 & 98.3 & 96.8 & 95.4 & 95.8 & 99.2 & 99.2 & 99.5 & 99.2 & 98.7 & 99.6 & 99.7 & 99.9 & 99.9 & 99.8 \\
	Ca2018& 99.4 & 100.0 & 99.4 & 99.0 & 98.6 & 99.6 & 100.0 & 99.7 & 99.3 & 98.8 & 99.8 & 100.0 & 99.8 & 99.4 & 98.9 \\
	Ho2018& 57.8 & 59.0 & 45.7 & 33.2 & 17.5 & 0.6 & 10.8 & 0.7 & 1.2 & 1.1 & 18.7 & 33.6 & 9.4 & 6.5 & 3.8 \\
	We2024& 99.7 & 99.7 & 99.7 & 99.7 & 99.6 & 99.8 & 99.8 & 99.8 & 99.8 & 99.8 & 99.9 & 99.9 & 99.9 & 99.9 & 99.9 \\
	Ma2015& 100.0 & 100.0 & 100.0 & 100.0 & 100.0 & 100.0 & 100.0 & 100.0 & 100.0 & 100.0 & 100.0 & 100.0 & 100.0 & 100.0 & 100.0 \\
	Mu2018& 100.0 & 100.0 & 100.0 & 100.0 & 100.0 & 100.0 & 100.0 & 100.0 & 100.0 & 100.0 & 100.0 & 100.0 & 100.0 & 100.0 & 100.0 \\
	\addlinespace[.2cm]
	\multicolumn{16}{l}{\quad\textit{Panel B: Specificity for $\mu_2$ }} \\
	Ours (full)& 98.6 & 98.6 & 98.6 & 98.2 & 98.1 & 99.3 & 99.3 & 99.6 & 99.5 & 99.3 & 99.7 & 99.7 & 99.8 & 99.8 & 99.7 \\
	Ours (partial)& 98.3 & 98.3 & 98.5 & 97.5 & 97.3 & 99.3 & 99.2 & 99.6 & 99.5 & 99.1 & 99.7 & 99.6 & 99.9 & 99.9 & 99.8 \\
	Ca2018& 99.4 & 100.0 & 99.4 & 99.0 & 98.6 & 99.6 & 100.0 & 99.6 & 99.3 & 98.8 & 99.8 & 100.0 & 99.8 & 99.4 & 99.0 \\
	Ho2018& 57.0 & 60.6 & 46.6 & 34.2 & 16.5 & 0.6 & 5.3 & 1.0 & 1.7 & 1.0 & 18.3 & 29.7 & 11.9 & 8.7 & 3.5 \\
	We2024& 99.7 & 99.7 & 99.7 & 99.7 & 99.6 & 99.8 & 99.8 & 99.8 & 99.8 & 99.8 & 99.9 & 99.9 & 99.9 & 99.9 & 99.9 \\
	Ma2015& 100.0 & 100.0 & 100.0 & 100.0 & 100.0 & 100.0 & 100.0 & 100.0 & 100.0 & 100.0 & 100.0 & 100.0 & 100.0 & 100.0 & 100.0 \\
	Mu2018& 100.0 & 100.0 & 100.0 & 100.0 & 100.0 & 100.0 & 100.0 & 100.0 & 100.0 & 100.0 & 100.0 & 100.0 & 100.0 & 100.0 & 100.0 \\
	\addlinespace[.2cm]
	\multicolumn{16}{l}{\quad\textit{Panel C: Specificity for $\mu_3$ }} \\
	Ours (full)& 97.5 & 97.6 & 98.1 & 97.8 & 97.8 & 98.9 & 98.8 & 99.4 & 99.3 & 99.3 & 99.5 & 99.5 & 99.7 & 99.8 & 99.7 \\
	Ours (partial)& 97.0 & 96.7 & 97.9 & 97.4 & 97.2 & 98.7 & 98.6 & 99.3 & 99.3 & 99.0 & 99.5 & 99.4 & 99.7 & 99.8 & 99.8 \\
	Ca2018& 99.6 & 99.6 & 98.7 & 98.5 & 98.4 & 99.7 & 99.7 & 98.9 & 98.9 & 98.7 & 99.7 & 99.8 & 99.3 & 99.3 & 98.9 \\
	Ho2018& 62.7 & 68.7 & 56.7 & 49.1 & 29.1 & 10.4 & 52.9 & 3.4 & 3.3 & 2.6 & 32.8 & 59.6 & 19.9 & 15.8 & 10.2 \\
	We2024& 99.7 & 99.7 & 99.7 & 99.7 & 99.6 & 99.8 & 99.8 & 99.8 & 99.8 & 99.8 & 99.9 & 99.9 & 99.9 & 99.9 & 99.9 \\
	Ma2015& 100.0 & 100.0 & 100.0 & 100.0 & 100.0 & 100.0 & 100.0 & 100.0 & 100.0 & 100.0 & 100.0 & 100.0 & 100.0 & 100.0 & 100.0 \\
	Mu2018& 100.0 & 100.0 & 100.0 & 100.0 & 100.0 & 100.0 & 100.0 & 100.0 & 100.0 & 100.0 & 100.0 & 100.0 & 100.0 & 100.0 & 100.0 \\
	\addlinespace[.2cm]
	\multicolumn{16}{l}{\quad\textit{Panel D: Specificity for $\mu_0$ }} \\
	Ours (full)& 98.3 & 98.6 & 98.8 & 98.5 & 97.9 & 99.3 & 99.3 & 99.6 & 99.5 & 99.4 & 99.7 & 99.7 & 99.8 & 99.8 & 99.8 \\
	Ours (partial)& 98.1 & 98.3 & 98.3 & 97.5 & 97.3 & 99.2 & 99.2 & 99.6 & 99.6 & 99.2 & 99.7 & 99.7 & 99.8 & 99.9 & 99.8 \\
	Ca2018& 99.5 & 100.0 & 99.3 & 99.0 & 98.6 & 99.7 & 100.0 & 99.7 & 99.3 & 98.8 & 99.8 & 100.0 & 99.8 & 99.4 & 99.0 \\
	Ho2018& 57.0 & 59.4 & 50.6 & 40.6 & 21.6 & 0.5 & 1.8 & 1.2 & 2.3 & 1.5 & 17.9 & 21.1 & 11.6 & 11.4 & 5.1 \\
	We2024& 99.7 & 99.7 & 99.7 & 99.6 & 99.6 & 99.8 & 99.8 & 99.8 & 99.8 & 99.8 & 99.9 & 99.9 & 99.9 & 99.9 & 99.9 \\
	Ma2015& 100.0 & 100.0 & 100.0 & 100.0 & 100.0 & 100.0 & 100.0 & 100.0 & 100.0 & 100.0 & 100.0 & 100.0 & 100.0 & 100.0 & 100.0 \\
	Mu2018& 100.0 & 100.0 & 100.0 & 100.0 & 100.0 & 100.0 & 100.0 & 100.0 & 100.0 & 100.0 & 100.0 & 100.0 & 100.0 & 100.0 & 100.0 \\
	\hline \hline 
 \end{tabular}
\end{footnotesize}
\end{table}
\begin{table}
\caption{Empirical specificity of various methods for AR errors, under different mean functions and error distributions.}
\begin{footnotesize}
\begin{tabular}{l| rrrrr | rrrrr | rrrrr}
	\hline \hline
	$n$ & \multicolumn{5}{c|}{50} & \multicolumn{5}{c|}{100} & \multicolumn{5}{c}{200} \\
	 $\mu$& $\Nc$ & $\Uc$ & $Exp$ & $Par_1$ & $Par_2$& $\Nc$ & $\Uc$ & $Exp$ & $Par_1$ & $Par_2$& $\Nc$ & $\Uc$ & $Exp$ & $Par_1$ & $Par_2$\\
	\hline
	\addlinespace[.2cm]
	\multicolumn{16}{l}{\quad\textit{Panel A: Specificity for $\mu_1$ }} \\
	Ours (full)& 98.4 & 98.3 & 98.6 & 98.3 & 98.0 & 99.3 & 99.3 & 99.5 & 99.4 & 99.2 & 99.7 & 99.7 & 99.8 & 99.8 & 99.7 \\
	Ours (partial)& 76.8 & 80.9 & 80.0 & 73.2 & 74.1 & 86.7 & 84.8 & 85.4 & 85.1 & 84.6 & 99.7 & 99.6 & 99.4 & 96.2 & 93.1 \\
	Ca2018& 99.5 & 99.8 & 99.1 & 98.9 & 98.4 & 99.7 & 99.8 & 99.4 & 99.2 & 98.6 & 99.8 & 99.9 & 99.6 & 99.3 & 98.7 \\
	Ho2018& 57.2 & 58.7 & 53.1 & 51.4 & 46.6 & 0.6 & 1.7 & 0.4 & 2.9 & 5.4 & 18.0 & 25.1 & 10.1 & 12.8 & 11.1 \\
	We2024& 99.7 & 99.7 & 99.7 & 99.7 & 99.6 & 99.8 & 99.8 & 99.8 & 99.8 & 99.8 & 99.9 & 99.9 & 99.9 & 99.9 & 99.9 \\
	Ma2015& 100.0 & 100.0 & 100.0 & 100.0 & 100.0 & 100.0 & 100.0 & 100.0 & 100.0 & 100.0 & 100.0 & 100.0 & 100.0 & 100.0 & 100.0 \\
	Mu2018& 100.0 & 100.0 & 100.0 & 100.0 & 100.0 & 100.0 & 100.0 & 100.0 & 100.0 & 100.0 & 100.0 & 100.0 & 100.0 & 100.0 & 100.0 \\
	\addlinespace[.2cm]
	\multicolumn{16}{l}{\quad\textit{Panel B: Specificity for $\mu_2$ }} \\
	Ours (full)& 98.3 & 98.2 & 98.8 & 98.4 & 97.7 & 99.3 & 99.3 & 99.4 & 99.3 & 99.2 & 99.7 & 99.7 & 99.8 & 99.8 & 99.7 \\
	Ours (partial)& 80.2 & 84.9 & 82.9 & 80.5 & 76.8 & 86.0 & 87.4 & 89.5 & 87.2 & 84.3 & 99.6 & 99.6 & 99.1 & 97.5 & 93.9 \\
	Ca2018& 99.4 & 99.8 & 99.1 & 98.9 & 98.5 & 99.7 & 99.8 & 99.4 & 99.2 & 98.6 & 99.8 & 99.9 & 99.5 & 99.3 & 98.7 \\
	Ho2018& 57.5 & 58.7 & 54.1 & 51.5 & 46.2 & 0.6 & 1.3 & 0.4 & 2.7 & 5.6 & 18.5 & 25.4 & 10.3 & 13.0 & 11.4 \\
	We2024& 99.7 & 99.7 & 99.7 & 99.7 & 99.6 & 99.8 & 99.8 & 99.9 & 99.8 & 99.8 & 99.9 & 99.9 & 99.9 & 99.9 & 99.9 \\
	Ma2015& 100.0 & 100.0 & 100.0 & 100.0 & 100.0 & 100.0 & 100.0 & 100.0 & 100.0 & 100.0 & 100.0 & 100.0 & 100.0 & 100.0 & 100.0 \\
	Mu2018& 100.0 & 100.0 & 100.0 & 100.0 & 100.0 & 100.0 & 100.0 & 100.0 & 100.0 & 100.0 & 100.0 & 100.0 & 100.0 & 100.0 & 100.0 \\
	\addlinespace[.2cm]
	\multicolumn{16}{l}{\quad\textit{Panel C: Specificity for $\mu_3$ }} \\
	Ours (full)& 97.6 & 97.7 & 98.2 & 97.9 & 97.7 & 98.9 & 98.9 & 99.2 & 99.2 & 99.1 & 99.5 & 99.5 & 99.7 & 99.7 & 99.6 \\
	Ours (partial)& 68.3 & 69.6 & 75.9 & 78.0 & 76.2 & 73.3 & 74.2 & 80.9 & 84.0 & 86.6 & 90.7 & 83.3 & 97.9 & 96.2 & 93.7 \\
	Ca2018& 99.6 & 99.6 & 99.5 & 98.9 & 98.4 & 99.8 & 99.8 & 99.7 & 99.1 & 98.6 & 99.9 & 99.9 & 99.7 & 99.2 & 98.7 \\
	Ho2018& 64.5 & 67.2 & 58.8 & 54.4 & 44.7 & 12.0 & 31.7 & 1.3 & 2.6 & 4.0 & 32.8 & 45.4 & 16.3 & 15.5 & 11.4 \\
	We2024& 99.6 & 99.5 & 99.7 & 99.6 & 99.6 & 99.8 & 99.8 & 99.8 & 99.8 & 99.8 & 99.9 & 99.9 & 99.9 & 99.9 & 99.9 \\
	Ma2015& 100.0 & 100.0 & 100.0 & 100.0 & 100.0 & 100.0 & 100.0 & 100.0 & 100.0 & 100.0 & 100.0 & 100.0 & 100.0 & 100.0 & 100.0 \\
	Mu2018& 100.0 & 100.0 & 100.0 & 100.0 & 100.0 & 100.0 & 100.0 & 100.0 & 100.0 & 100.0 & 100.0 & 100.0 & 100.0 & 100.0 & 100.0 \\
	\addlinespace[.2cm]
	\multicolumn{16}{l}{\quad\textit{Panel D: Specificity for $\mu_0$ }} \\
	Ours (full)& 98.3 & 98.0 & 98.6 & 98.3 & 97.9 & 99.3 & 99.3 & 99.5 & 99.5 & 99.1 & 99.7 & 99.7 & 99.8 & 99.8 & 99.7 \\
	Ours (partial)& 79.7 & 83.0 & 83.5 & 80.7 & 77.0 & 86.0 & 86.8 & 89.2 & 86.1 & 84.5 & 99.6 & 99.6 & 99.6 & 97.2 & 94.7 \\
	Ca2018& 99.5 & 99.8 & 99.1 & 98.9 & 98.4 & 99.7 & 99.8 & 99.4 & 99.2 & 98.6 & 99.8 & 99.9 & 99.5 & 99.3 & 98.7 \\
	Ho2018& 57.7 & 59.2 & 53.5 & 51.0 & 46.5 & 0.7 & 1.2 & 0.5 & 2.7 & 5.4 & 18.6 & 25.1 & 9.5 & 13.0 & 11.4 \\
	We2024& 99.7 & 99.7 & 99.7 & 99.7 & 99.6 & 99.8 & 99.8 & 99.8 & 99.8 & 99.8 & 99.9 & 99.9 & 99.9 & 99.9 & 99.9 \\
	Ma2015& 100.0 & 100.0 & 100.0 & 100.0 & 100.0 & 100.0 & 100.0 & 100.0 & 100.0 & 100.0 & 100.0 & 100.0 & 100.0 & 100.0 & 100.0 \\
	Mu2018& 100.0 & 100.0 & 100.0 & 100.0 & 100.0 & 100.0 & 100.0 & 100.0 & 100.0 & 100.0 & 100.0 & 100.0 & 100.0 & 100.0 & 100.0 \\
	\hline \hline 
 \end{tabular}
\end{footnotesize}
\end{table}
\begin{table}
\caption{Empirical specificity of various methods for MA errors, under different mean functions and error distributions.}
\begin{footnotesize}
\begin{tabular}{l| rrrrr | rrrrr | rrrrr}
	\hline \hline
	$n$ & \multicolumn{5}{c|}{50} & \multicolumn{5}{c|}{100} & \multicolumn{5}{c}{200} \\
	 $\mu$& $\Nc$ & $\Uc$ & $Exp$ & $Par_1$ & $Par_2$& $\Nc$ & $\Uc$ & $Exp$ & $Par_1$ & $Par_2$& $\Nc$ & $\Uc$ & $Exp$ & $Par_1$ & $Par_2$\\
	\hline
	\addlinespace[.2cm]
	\multicolumn{16}{l}{\quad\textit{Panel A: Specificity for $\mu_1$ }} \\
	Ours (full)& 98.2 & 98.3 & 98.5 & 98.2 & 97.7 & 99.3 & 99.2 & 99.4 & 99.1 & 99.0 & 99.7 & 99.7 & 99.8 & 99.7 & 99.7 \\
	Ours (partial)& 86.3 & 87.9 & 89.2 & 82.1 & 82.2 & 95.6 & 95.5 & 96.8 & 94.3 & 95.0 & 99.7 & 99.6 & 99.8 & 99.4 & 98.5 \\
	Ca2018& 99.5 & 99.7 & 99.2 & 98.9 & 98.4 & 99.7 & 99.8 & 99.5 & 99.1 & 98.6 & 99.8 & 99.9 & 99.7 & 99.3 & 98.7 \\
	Ho2018& 57.3 & 58.4 & 51.5 & 46.6 & 38.9 & 0.6 & 1.2 & 0.6 & 2.7 & 4.5 & 18.5 & 22.3 & 11.5 & 11.3 & 10.9 \\
	We2024& 99.7 & 99.7 & 99.7 & 99.6 & 99.6 & 99.8 & 99.8 & 99.8 & 99.8 & 99.8 & 99.9 & 99.9 & 99.9 & 99.9 & 99.9 \\
	Ma2015& 100.0 & 100.0 & 100.0 & 100.0 & 100.0 & 100.0 & 100.0 & 100.0 & 100.0 & 100.0 & 100.0 & 100.0 & 100.0 & 100.0 & 100.0 \\
	Mu2018& 100.0 & 100.0 & 100.0 & 100.0 & 100.0 & 100.0 & 100.0 & 100.0 & 100.0 & 100.0 & 100.0 & 100.0 & 100.0 & 100.0 & 100.0 \\
	\addlinespace[.2cm]
	\multicolumn{16}{l}{\quad\textit{Panel B: Specificity for $\mu_2$ }} \\
	Ours (full)& 98.2 & 98.3 & 98.4 & 98.1 & 97.9 & 99.2 & 99.3 & 99.4 & 99.3 & 99.1 & 99.7 & 99.7 & 99.8 & 99.7 & 99.6 \\
	Ours (partial)& 91.3 & 91.6 & 91.6 & 90.3 & 88.8 & 97.1 & 97.2 & 98.0 & 97.2 & 95.2 & 99.6 & 99.6 & 99.8 & 99.7 & 98.4 \\
	Ca2018& 99.5 & 99.7 & 99.2 & 98.9 & 98.4 & 99.7 & 99.8 & 99.4 & 99.1 & 98.6 & 99.8 & 99.9 & 99.7 & 99.3 & 98.7 \\
	Ho2018& 57.1 & 58.7 & 51.7 & 46.9 & 37.8 & 0.5 & 1.4 & 0.8 & 2.9 & 4.0 & 17.7 & 22.4 & 11.3 & 11.2 & 11.3 \\
	We2024& 99.7 & 99.7 & 99.7 & 99.6 & 99.6 & 99.8 & 99.8 & 99.8 & 99.8 & 99.8 & 99.9 & 99.9 & 99.9 & 99.9 & 99.9 \\
	Ma2015& 100.0 & 100.0 & 100.0 & 100.0 & 100.0 & 100.0 & 100.0 & 100.0 & 100.0 & 100.0 & 100.0 & 100.0 & 100.0 & 100.0 & 100.0 \\
	Mu2018& 100.0 & 100.0 & 100.0 & 100.0 & 100.0 & 100.0 & 100.0 & 100.0 & 100.0 & 100.0 & 100.0 & 100.0 & 100.0 & 100.0 & 100.0 \\
	\addlinespace[.2cm]
	\multicolumn{16}{l}{\quad\textit{Panel C: Specificity for $\mu_3$ }} \\
	Ours (full)& 97.5 & 97.5 & 98.0 & 97.7 & 97.3 & 99.0 & 98.8 & 99.2 & 99.1 & 99.1 & 99.5 & 99.5 & 99.7 & 99.7 & 99.6 \\
	Ours (partial)& 87.8 & 87.1 & 89.9 & 87.4 & 89.1 & 95.8 & 94.9 & 96.4 & 95.8 & 96.0 & 99.3 & 98.9 & 99.6 & 99.5 & 98.7 \\
	Ca2018& 99.6 & 99.7 & 99.4 & 98.9 & 98.4 & 99.8 & 99.8 & 99.6 & 99.1 & 98.5 & 99.9 & 99.9 & 99.6 & 99.2 & 98.7 \\
	Ho2018& 63.2 & 65.5 & 57.0 & 51.2 & 37.7 & 7.5 & 25.3 & 1.6 & 2.4 & 3.5 & 29.5 & 43.3 & 17.9 & 12.9 & 9.9 \\
	We2024& 99.6 & 99.6 & 99.7 & 99.6 & 99.6 & 99.8 & 99.8 & 99.8 & 99.8 & 99.8 & 99.9 & 99.9 & 99.9 & 99.9 & 99.9 \\
	Ma2015& 100.0 & 100.0 & 100.0 & 100.0 & 100.0 & 100.0 & 100.0 & 100.0 & 100.0 & 100.0 & 100.0 & 100.0 & 100.0 & 100.0 & 100.0 \\
	Mu2018& 100.0 & 100.0 & 100.0 & 100.0 & 100.0 & 100.0 & 100.0 & 100.0 & 100.0 & 100.0 & 100.0 & 100.0 & 100.0 & 100.0 & 100.0 \\
	\addlinespace[.2cm]
	\multicolumn{16}{l}{\quad\textit{Panel D: Specificity for $\mu_0$ }} \\
	Ours (full)& 98.2 & 98.1 & 98.6 & 98.0 & 97.7 & 99.3 & 99.2 & 99.4 & 99.3 & 99.0 & 99.7 & 99.7 & 99.8 & 99.8 & 99.7 \\
	Ours (partial)& 92.2 & 92.1 & 92.8 & 88.9 & 87.0 & 96.9 & 97.7 & 97.5 & 96.6 & 94.4 & 99.6 & 99.6 & 99.7 & 99.4 & 99.4 \\
	Ca2018& 99.4 & 99.7 & 99.2 & 98.9 & 98.5 & 99.6 & 99.8 & 99.5 & 99.1 & 98.6 & 99.8 & 99.9 & 99.7 & 99.3 & 98.8 \\
	Ho2018& 57.3 & 58.4 & 51.9 & 46.0 & 38.6 & 0.6 & 1.2 & 0.6 & 2.3 & 4.0 & 18.2 & 22.5 & 11.4 & 11.1 & 10.8 \\
	We2024& 99.7 & 99.7 & 99.7 & 99.6 & 99.6 & 99.8 & 99.8 & 99.9 & 99.8 & 99.8 & 99.9 & 99.9 & 99.9 & 99.9 & 99.9 \\
	Ma2015& 100.0 & 100.0 & 100.0 & 100.0 & 100.0 & 100.0 & 100.0 & 100.0 & 100.0 & 100.0 & 100.0 & 100.0 & 100.0 & 100.0 & 100.0 \\
	Mu2018& 100.0 & 100.0 & 100.0 & 100.0 & 100.0 & 100.0 & 100.0 & 100.0 & 100.0 & 100.0 & 100.0 & 100.0 & 100.0 & 100.0 & 100.0 \\
	\hline \hline 
 \end{tabular}
\end{footnotesize}
\end{table}
\begin{table}
\caption{Empirical specificity and sensitivity of various methods for IID errors, under different mean functions and error distributions.}
\begin{footnotesize}
\begin{tabular}{l| rrrrr | rrrrr | rrrrr}
	\hline \hline
	$n$ & \multicolumn{5}{c|}{50} & \multicolumn{5}{c|}{100} & \multicolumn{5}{c}{200} \\
	 $\mu$& $\Nc$ & $\Uc$ & $Exp$ & $Par_1$ & $Par_2$& $\Nc$ & $\Uc$ & $Exp$ & $Par_1$ & $Par_2$& $\Nc$ & $\Uc$ & $Exp$ & $Par_1$ & $Par_2$\\
	\hline
	\addlinespace[.2cm]
	\multicolumn{16}{l}{\quad\textit{Panel A: Specificity for $\mu_1$ }} \\
	Ours (full)& 86.8 & 95.2 & 85.5 & 89.4 & 88.3 & 87.4 & 96.7 & 87.6 & 93.4 & 92.6 & 87.2 & 97.3 & 88.0 & 95.7 & 94.2 \\
	Ours (partial)& 90.2 & 94.9 & 89.5 & 87.3 & 91.3 & 98.3 & 97.9 & 98.3 & 98.0 & 98.7 & 99.6 & 98.8 & 99.7 & 99.6 & 99.1 \\
	Ca2018& 99.9 & 100.0 & 99.9 & 99.7 & 99.7 & 100.0 & 100.0 & 99.9 & 99.8 & 99.8 & 100.0 & 100.0 & 99.9 & 99.8 & 99.9 \\
	Ho2018& 72.5 & 72.2 & 70.5 & 57.2 & 60.8 & 67.9 & 66.6 & 66.5 & 43.6 & 62.3 & 71.3 & 71.5 & 71.2 & 61.8 & 69.1 \\
	We2024& 99.7 & 99.7 & 99.7 & 99.7 & 99.7 & 99.9 & 99.8 & 99.8 & 99.8 & 99.8 & 99.9 & 99.9 & 99.9 & 99.9 & 99.9 \\
	Ma2015& 100.0 & 100.0 & 100.0 & 100.0 & 100.0 & 100.0 & 100.0 & 100.0 & 100.0 & 100.0 & 100.0 & 100.0 & 100.0 & 100.0 & 100.0 \\
	Mu2018& 100.0 & 100.0 & 100.0 & 100.0 & 100.0 & 100.0 & 100.0 & 100.0 & 100.0 & 100.0 & 100.0 & 100.0 & 100.0 & 100.0 & 100.0 \\
	\addlinespace[.2cm]
	\multicolumn{16}{l}{\quad\textit{Panel B: Specificity for $\mu_2$ }} \\
	Ours (full)& 87.2 & 95.7 & 87.8 & 92.5 & 90.2 & 87.5 & 96.8 & 88.5 & 94.9 & 92.9 & 87.3 & 97.4 & 88.1 & 95.9 & 94.3 \\
	Ours (partial)& 95.5 & 96.3 & 95.6 & 94.9 & 94.7 & 99.0 & 98.2 & 99.2 & 98.9 & 98.3 & 99.7 & 98.9 & 99.5 & 99.5 & 99.2 \\
	Ca2018& 99.9 & 100.0 & 99.9 & 99.7 & 99.7 & 100.0 & 100.0 & 99.9 & 99.8 & 99.8 & 100.0 & 100.0 & 99.9 & 99.9 & 99.9 \\
	Ho2018& 73.0 & 74.8 & 71.6 & 56.7 & 57.5 & 68.8 & 70.8 & 67.9 & 41.0 & 58.9 & 71.8 & 73.7 & 72.3 & 62.0 & 68.2 \\
	We2024& 99.7 & 99.7 & 99.7 & 99.7 & 99.7 & 99.8 & 99.8 & 99.8 & 99.8 & 99.8 & 99.9 & 99.9 & 99.9 & 99.9 & 99.9 \\
	Ma2015& 100.0 & 100.0 & 100.0 & 100.0 & 100.0 & 100.0 & 100.0 & 100.0 & 100.0 & 100.0 & 100.0 & 100.0 & 100.0 & 100.0 & 100.0 \\
	Mu2018& 100.0 & 100.0 & 100.0 & 100.0 & 100.0 & 100.0 & 100.0 & 100.0 & 100.0 & 100.0 & 100.0 & 100.0 & 100.0 & 100.0 & 100.0 \\
	\addlinespace[.2cm]
	\multicolumn{16}{l}{\quad\textit{Panel C: Specificity for $\mu_3$ }} \\
	Ours (full)& 86.8 & 94.7 & 87.3 & 91.9 & 90.1 & 87.1 & 96.4 & 88.4 & 94.8 & 92.8 & 87.0 & 97.1 & 87.8 & 96.0 & 93.8 \\
	Ours (partial)& 95.0 & 95.2 & 95.0 & 95.3 & 95.4 & 98.6 & 97.7 & 98.3 & 98.9 & 98.3 & 99.4 & 98.5 & 99.4 & 99.5 & 99.3 \\
	Ca2018& 99.9 & 99.7 & 99.8 & 99.7 & 99.7 & 99.9 & 99.9 & 99.9 & 99.8 & 99.8 & 100.0 & 99.9 & 99.9 & 99.8 & 99.9 \\
	Ho2018& 67.8 & 70.6 & 63.6 & 42.5 & 53.3 & 63.7 & 65.8 & 60.8 & 23.5 & 57.1 & 69.1 & 70.8 & 68.0 & 51.5 & 65.0 \\
	We2024& 99.7 & 99.7 & 99.7 & 99.7 & 99.7 & 99.9 & 99.9 & 99.8 & 99.8 & 99.8 & 99.9 & 99.9 & 99.9 & 99.9 & 99.9 \\
	Ma2015& 100.0 & 100.0 & 100.0 & 100.0 & 100.0 & 100.0 & 100.0 & 100.0 & 100.0 & 100.0 & 100.0 & 100.0 & 100.0 & 100.0 & 100.0 \\
	Mu2018& 100.0 & 100.0 & 100.0 & 100.0 & 100.0 & 100.0 & 100.0 & 100.0 & 100.0 & 100.0 & 100.0 & 100.0 & 100.0 & 100.0 & 100.0 \\
	\addlinespace[.2cm]
	\multicolumn{16}{l}{\quad\textit{Panel D: Specificity for $\mu_0$ }} \\
	Ours (full)& 87.6 & 96.0 & 87.4 & 92.3 & 90.6 & 87.7 & 96.8 & 88.5 & 94.9 & 92.9 & 87.4 & 97.3 & 88.1 & 95.5 & 93.7 \\
	Ours (partial)& 95.4 & 96.7 & 96.1 & 96.0 & 95.4 & 99.2 & 98.2 & 98.7 & 98.9 & 98.7 & 99.7 & 98.8 & 99.4 & 99.6 & 99.3 \\
	Ca2018& 99.9 & 100.0 & 99.9 & 99.7 & 99.7 & 100.0 & 100.0 & 99.9 & 99.8 & 99.8 & 100.0 & 100.0 & 99.9 & 99.9 & 99.9 \\
	Ho2018& 73.3 & 75.8 & 72.3 & 58.0 & 58.0 & 69.1 & 72.3 & 68.5 & 41.1 & 59.1 & 71.8 & 74.4 & 72.8 & 62.9 & 68.5 \\
	We2024& 99.7 & 99.7 & 99.7 & 99.7 & 99.7 & 99.8 & 99.9 & 99.9 & 99.8 & 99.8 & 99.9 & 99.9 & 99.9 & 99.9 & 99.9 \\
	Ma2015& 100.0 & 100.0 & 100.0 & 100.0 & 100.0 & 100.0 & 100.0 & 100.0 & 100.0 & 100.0 & 100.0 & 100.0 & 100.0 & 100.0 & 100.0 \\
	Mu2018& 100.0 & 100.0 & 100.0 & 100.0 & 100.0 & 100.0 & 100.0 & 100.0 & 100.0 & 100.0 & 100.0 & 100.0 & 100.0 & 100.0 & 100.0 \\
	\addlinespace[.2cm]
	\multicolumn{16}{l}{\quad\textit{Panel E: Sensitivity for $\mu_1$ }} \\
	Ours (full)& 80.4 & 32.7 & 79.1 & 52.4 & 63.8 & 93.8 & 48.7 & 91.3 & 60.9 & 61.2 & 97.5 & 61.1 & 95.4 & 60.6 & 58.0 \\
	Ours (partial)& 82.9 & 31.6 & 81.4 & 56.5 & 65.5 & 94.9 & 47.0 & 92.6 & 61.2 & 62.0 & 99.1 & 59.8 & 96.8 & 61.5 & 58.9 \\
	Ca2018& 82.9 & 41.5 & 82.9 & 78.3 & 79.4 & 86.1 & 38.5 & 85.7 & 82.6 & 84.0 & 87.2 & 40.3 & 87.1 & 84.0 & 87.2 \\
	Ho2018& 100.0 & 99.3 & 100.0 & 100.0 & 100.0 & 100.0 & 99.9 & 100.0 & 100.0 & 94.0 & 100.0 & 100.0 & 99.6 & 100.0 & 100.0 \\
	We2024& 2.0 & 0.5 & 3.1 & 2.4 & 3.1 & 1.1 & 0.2 & 1.6 & 1.3 & 1.7 & 0.6 & 0.1 & 0.9 & 0.7 & 0.9 \\
	Ma2015& 0.0 & 0.0 & 0.0 & 0.0 & 0.0 & 0.0 & 0.0 & 0.0 & 0.0 & 0.0 & 0.0 & 0.0 & 0.0 & 0.0 & 0.0 \\
	Mu2018& 0.0 & 0.0 & 0.0 & 0.0 & 0.2 & 0.0 & 0.0 & 0.0 & 0.0 & 0.0 & 0.0 & 0.0 & 0.0 & 0.0 & 0.0 \\
	\addlinespace[.2cm]
	\multicolumn{16}{l}{\quad\textit{Panel F: Sensitivity for $\mu_2$ }} \\
	Ours (full)& 90.6 & 49.3 & 87.0 & 64.6 & 63.9 & 96.4 & 56.1 & 92.6 & 59.6 & 61.1 & 97.6 & 60.6 & 95.2 & 60.1 & 58.2 \\
	Ours (partial)& 92.9 & 47.9 & 89.0 & 66.5 & 65.9 & 97.7 & 54.3 & 94.0 & 60.5 & 62.5 & 99.1 & 59.3 & 96.8 & 61.0 & 59.1 \\
	Ca2018& 84.1 & 45.4 & 84.0 & 79.6 & 80.9 & 86.4 & 46.7 & 86.3 & 82.6 & 85.0 & 87.4 & 46.5 & 87.4 & 84.7 & 87.8 \\
	Ho2018& 100.0 & 99.4 & 100.0 & 99.9 & 99.9 & 100.0 & 100.0 & 100.0 & 100.0 & 91.4 & 100.0 & 100.0 & 98.9 & 100.0 & 100.0 \\
	We2024& 2.1 & 0.5 & 3.1 & 2.4 & 3.0 & 1.1 & 0.2 & 1.7 & 1.3 & 1.7 & 0.6 & 0.1 & 0.8 & 0.7 & 0.9 \\
	Ma2015& 0.0 & 0.0 & 0.0 & 0.0 & 1.0 & 0.0 & 0.0 & 0.0 & 0.0 & 0.0 & 0.0 & 0.0 & 0.0 & 0.0 & 0.0 \\
	Mu2018& 0.0 & 0.0 & 0.0 & 0.0 & 0.3 & 0.0 & 0.0 & 0.0 & 0.0 & 0.0 & 0.0 & 0.0 & 0.0 & 0.0 & 0.0 \\
	\addlinespace[.2cm]
	\multicolumn{16}{l}{\quad\textit{Panel G: Sensitivity for $\mu_3$ }} \\
	Ours (full)& 90.7 & 49.0 & 86.5 & 62.9 & 64.5 & 96.5 & 56.3 & 91.7 & 60.9 & 61.7 & 97.5 & 60.6 & 95.5 & 58.4 & 60.4 \\
	Ours (partial)& 91.6 & 47.1 & 88.1 & 64.1 & 66.5 & 97.6 & 54.6 & 93.5 & 61.5 & 62.9 & 99.1 & 59.2 & 97.2 & 59.2 & 61.2 \\
	Ca2018& 77.3 & 33.2 & 81.3 & 74.1 & 80.0 & 83.5 & 30.6 & 85.1 & 80.8 & 84.3 & 86.0 & 27.9 & 86.9 & 83.6 & 87.3 \\
	Ho2018& 100.0 & 98.7 & 100.0 & 100.0 & 100.0 & 100.0 & 99.7 & 100.0 & 100.0 & 93.6 & 99.9 & 99.9 & 99.8 & 100.0 & 100.0 \\
	We2024& 2.1 & 0.5 & 3.1 & 2.4 & 3.1 & 1.1 & 0.2 & 1.6 & 1.3 & 1.7 & 0.6 & 0.1 & 0.9 & 0.7 & 0.9 \\
	Ma2015& 0.0 & 0.0 & 0.0 & 0.1 & 0.9 & 0.0 & 0.0 & 0.0 & 0.0 & 0.1 & 0.0 & 0.0 & 0.0 & 0.0 & 0.0 \\
	Mu2018& 0.0 & 0.0 & 0.0 & 0.0 & 0.3 & 0.0 & 0.0 & 0.0 & 0.0 & 0.0 & 0.0 & 0.0 & 0.0 & 0.0 & 0.0 \\
	\addlinespace[.2cm]
	\multicolumn{16}{l}{\quad\textit{Panel H: Sensitivity for $\mu_0$ }} \\
	Ours (full)& 90.8 & 48.7 & 87.7 & 63.2 & 63.9 & 96.1 & 55.7 & 91.5 & 61.0 & 59.5 & 97.5 & 61.5 & 94.6 & 61.2 & 59.7 \\
	Ours (partial)& 92.9 & 46.9 & 89.3 & 64.7 & 65.6 & 97.6 & 54.3 & 93.2 & 62.1 & 60.6 & 99.1 & 60.3 & 96.3 & 62.0 & 60.5 \\
	Ca2018& 84.6 & 46.4 & 84.4 & 80.1 & 80.7 & 86.9 & 50.6 & 86.5 & 82.6 & 85.1 & 87.5 & 52.6 & 87.7 & 84.8 & 87.7 \\
	Ho2018& 100.0 & 99.4 & 100.0 & 99.8 & 99.6 & 100.0 & 100.0 & 100.0 & 100.0 & 90.9 & 100.0 & 100.0 & 99.1 & 100.0 & 100.0 \\
	We2024& 2.1 & 0.4 & 3.1 & 2.3 & 3.2 & 1.1 & 0.2 & 1.7 & 1.3 & 1.7 & 0.6 & 0.1 & 0.9 & 0.7 & 0.9 \\
	Ma2015& 0.0 & 0.0 & 0.0 & 0.0 & 0.7 & 0.0 & 0.0 & 0.0 & 0.0 & 0.0 & 0.0 & 0.0 & 0.0 & 0.0 & 0.0 \\
	Mu2018& 0.0 & 0.0 & 0.0 & 0.0 & 0.3 & 0.0 & 0.0 & 0.0 & 0.0 & 0.0 & 0.0 & 0.0 & 0.0 & 0.0 & 0.0 \\
	\hline \hline 
 \end{tabular}
\end{footnotesize}
\end{table}
\begin{table}
\caption{Empirical specificity and sensitivity of various methods for AR errors, under different mean functions and error distributions.}
\begin{footnotesize}
\begin{tabular}{l| rrrrr | rrrrr | rrrrr}
	\hline \hline
	$n$ & \multicolumn{5}{c|}{50} & \multicolumn{5}{c|}{100} & \multicolumn{5}{c}{200} \\
	 $\mu$& $\Nc$ & $\Uc$ & $Exp$ & $Par_1$ & $Par_2$& $\Nc$ & $\Uc$ & $Exp$ & $Par_1$ & $Par_2$& $\Nc$ & $\Uc$ & $Exp$ & $Par_1$ & $Par_2$\\
	\hline
	\addlinespace[.2cm]
	\multicolumn{16}{l}{\quad\textit{Panel A: Specificity for $\mu_1$ }} \\
	Ours (full)& 86.5 & 92.5 & 86.9 & 90.7 & 89.3 & 87.3 & 93.9 & 86.8 & 92.1 & 91.1 & 85.3 & 94.0 & 85.3 & 92.8 & 91.0 \\
	Ours (partial)& 34.3 & 44.3 & 39.5 & 42.9 & 50.1 & 44.3 & 47.2 & 44.6 & 60.6 & 59.1 & 58.1 & 58.1 & 58.4 & 83.6 & 72.9 \\
	Ca2018& 99.9 & 100.0 & 99.9 & 99.7 & 99.7 & 99.9 & 100.0 & 99.9 & 99.8 & 99.8 & 100.0 & 100.0 & 99.9 & 99.9 & 99.9 \\
	Ho2018& 72.4 & 68.4 & 71.5 & 62.2 & 59.3 & 68.3 & 51.8 & 68.6 & 54.8 & 63.0 & 71.5 & 61.0 & 72.3 & 65.7 & 69.5 \\
	We2024& 99.7 & 99.7 & 99.7 & 99.7 & 99.7 & 99.8 & 99.8 & 99.9 & 99.8 & 99.8 & 99.9 & 99.9 & 99.9 & 99.9 & 99.9 \\
	Ma2015& 100.0 & 100.0 & 100.0 & 100.0 & 100.0 & 100.0 & 100.0 & 100.0 & 100.0 & 100.0 & 100.0 & 100.0 & 100.0 & 100.0 & 100.0 \\
	Mu2018& 100.0 & 100.0 & 100.0 & 100.0 & 100.0 & 100.0 & 100.0 & 100.0 & 100.0 & 100.0 & 100.0 & 100.0 & 100.0 & 100.0 & 100.0 \\
	\addlinespace[.2cm]
	\multicolumn{16}{l}{\quad\textit{Panel B: Specificity for $\mu_2$ }} \\
	Ours (full)& 86.6 & 92.9 & 86.5 & 90.5 & 89.5 & 86.9 & 94.6 & 87.2 & 92.3 & 91.1 & 85.4 & 94.0 & 85.3 & 93.0 & 91.5 \\
	Ours (partial)& 46.6 & 54.5 & 50.6 & 59.0 & 59.3 & 53.4 & 56.8 & 56.3 & 66.2 & 64.8 & 57.5 & 60.2 & 60.0 & 81.5 & 74.5 \\
	Ca2018& 99.9 & 100.0 & 99.9 & 99.8 & 99.7 & 99.9 & 100.0 & 99.9 & 99.8 & 99.8 & 100.0 & 100.0 & 99.9 & 99.9 & 99.9 \\
	Ho2018& 73.0 & 68.4 & 72.4 & 61.0 & 57.1 & 68.9 & 52.4 & 69.4 & 52.9 & 60.2 & 71.8 & 60.9 & 73.2 & 66.2 & 69.1 \\
	We2024& 99.7 & 99.7 & 99.7 & 99.7 & 99.7 & 99.8 & 99.8 & 99.8 & 99.8 & 99.8 & 99.9 & 99.9 & 99.9 & 99.9 & 99.9 \\
	Ma2015& 100.0 & 100.0 & 100.0 & 100.0 & 100.0 & 100.0 & 100.0 & 100.0 & 100.0 & 100.0 & 100.0 & 100.0 & 100.0 & 100.0 & 100.0 \\
	Mu2018& 100.0 & 100.0 & 100.0 & 100.0 & 100.0 & 100.0 & 100.0 & 100.0 & 100.0 & 100.0 & 100.0 & 100.0 & 100.0 & 100.0 & 100.0 \\
	\addlinespace[.2cm]
	\multicolumn{16}{l}{\quad\textit{Panel C: Specificity for $\mu_3$ }} \\
	Ours (full)& 86.5 & 92.2 & 86.7 & 89.9 & 88.7 & 86.9 & 93.9 & 86.6 & 92.2 & 91.2 & 84.7 & 94.4 & 85.7 & 92.8 & 91.5 \\
	Ours (partial)& 47.6 & 53.8 & 47.7 & 56.8 & 60.8 & 52.4 & 54.1 & 54.7 & 67.9 & 66.4 & 53.8 & 62.8 & 63.1 & 79.2 & 76.2 \\
	Ca2018& 99.9 & 99.8 & 99.8 & 99.7 & 99.7 & 99.9 & 99.9 & 99.9 & 99.8 & 99.8 & 100.0 & 100.0 & 99.9 & 99.9 & 99.9 \\
	Ho2018& 68.2 & 66.0 & 65.4 & 51.9 & 55.7 & 64.7 & 48.5 & 63.1 & 40.0 & 59.6 & 69.3 & 59.3 & 69.4 & 59.4 & 67.1 \\
	We2024& 99.7 & 99.7 & 99.7 & 99.7 & 99.7 & 99.8 & 99.9 & 99.9 & 99.8 & 99.8 & 99.9 & 99.9 & 99.9 & 99.9 & 99.9 \\
	Ma2015& 100.0 & 100.0 & 100.0 & 100.0 & 100.0 & 100.0 & 100.0 & 100.0 & 100.0 & 100.0 & 100.0 & 100.0 & 100.0 & 100.0 & 100.0 \\
	Mu2018& 100.0 & 100.0 & 100.0 & 100.0 & 100.0 & 100.0 & 100.0 & 100.0 & 100.0 & 100.0 & 100.0 & 100.0 & 100.0 & 100.0 & 100.0 \\
	\addlinespace[.2cm]
	\multicolumn{16}{l}{\quad\textit{Panel D: Specificity for $\mu_0$ }} \\
	Ours (full)& 86.2 & 92.8 & 87.0 & 90.4 & 88.6 & 87.4 & 94.6 & 87.4 & 92.3 & 90.8 & 85.3 & 94.3 & 85.3 & 93.0 & 91.5 \\
	Ours (partial)& 48.5 & 54.6 & 51.8 & 56.9 & 60.4 & 50.5 & 59.9 & 57.6 & 64.0 & 63.0 & 55.9 & 61.8 & 70.8 & 84.2 & 72.6 \\
	Ca2018& 99.9 & 100.0 & 99.9 & 99.8 & 99.7 & 100.0 & 100.0 & 99.9 & 99.8 & 99.8 & 100.0 & 100.0 & 99.9 & 99.9 & 99.9 \\
	Ho2018& 73.1 & 68.6 & 73.1 & 62.2 & 58.7 & 69.2 & 52.0 & 70.1 & 53.1 & 60.6 & 71.8 & 61.1 & 73.6 & 66.2 & 69.5 \\
	We2024& 99.7 & 99.7 & 99.7 & 99.7 & 99.7 & 99.8 & 99.8 & 99.8 & 99.8 & 99.8 & 99.9 & 99.9 & 99.9 & 99.9 & 99.9 \\
	Ma2015& 100.0 & 100.0 & 100.0 & 100.0 & 100.0 & 100.0 & 100.0 & 100.0 & 100.0 & 100.0 & 100.0 & 100.0 & 100.0 & 100.0 & 100.0 \\
	Mu2018& 100.0 & 100.0 & 100.0 & 100.0 & 100.0 & 100.0 & 100.0 & 100.0 & 100.0 & 100.0 & 100.0 & 100.0 & 100.0 & 100.0 & 100.0 \\
	\addlinespace[.2cm]
	\multicolumn{16}{l}{\quad\textit{Panel E: Sensitivity for $\mu_1$ }} \\
	Ours (full)& 76.4 & 38.1 & 75.4 & 52.0 & 56.4 & 82.3 & 42.1 & 77.9 & 50.9 & 53.3 & 77.3 & 31.6 & 74.1 & 38.3 & 47.2 \\
	Ours (partial)& 87.3 & 64.4 & 81.0 & 68.1 & 70.4 & 86.7 & 66.1 & 84.8 & 55.4 & 64.9 & 87.5 & 62.5 & 88.6 & 45.2 & 60.5 \\
	Ca2018& 82.2 & 42.5 & 83.7 & 79.4 & 80.4 & 85.8 & 34.8 & 86.1 & 83.0 & 84.5 & 87.2 & 24.6 & 87.3 & 84.7 & 87.3 \\
	Ho2018& 100.0 & 99.6 & 100.0 & 100.0 & 99.9 & 100.0 & 100.0 & 100.0 & 100.0 & 93.0 & 100.0 & 100.0 & 99.1 & 100.0 & 100.0 \\
	We2024& 3.2 & 0.6 & 3.7 & 3.0 & 3.4 & 1.6 & 0.2 & 1.9 & 1.6 & 1.8 & 0.8 & 0.1 & 0.9 & 0.8 & 0.9 \\
	Ma2015& 0.0 & 0.0 & 0.0 & 0.0 & 0.0 & 0.0 & 0.0 & 0.0 & 0.0 & 0.0 & 0.0 & 0.0 & 0.0 & 0.0 & 0.0 \\
	Mu2018& 0.0 & 0.0 & 0.0 & 0.0 & 0.1 & 0.0 & 0.0 & 0.0 & 0.0 & 0.0 & 0.0 & 0.0 & 0.0 & 0.0 & 0.0 \\
	\addlinespace[.2cm]
	\multicolumn{16}{l}{\quad\textit{Panel F: Sensitivity for $\mu_2$ }} \\
	Ours (full)& 82.2 & 42.3 & 79.0 & 59.2 & 55.9 & 85.2 & 45.5 & 82.8 & 54.1 & 55.6 & 77.3 & 32.8 & 74.2 & 42.7 & 44.1 \\
	Ours (partial)& 89.1 & 63.1 & 86.1 & 70.3 & 66.9 & 90.7 & 63.7 & 87.9 & 61.2 & 65.2 & 89.7 & 61.8 & 88.2 & 51.3 & 56.7 \\
	Ca2018& 83.6 & 43.6 & 84.5 & 79.6 & 80.2 & 86.1 & 36.5 & 86.7 & 83.2 & 85.4 & 87.5 & 25.0 & 88.0 & 85.3 & 87.8 \\
	Ho2018& 100.0 & 99.7 & 100.0 & 100.0 & 99.9 & 100.0 & 100.0 & 100.0 & 100.0 & 92.5 & 100.0 & 100.0 & 98.8 & 100.0 & 100.0 \\
	We2024& 3.3 & 0.6 & 3.6 & 2.9 & 3.3 & 1.6 & 0.3 & 1.8 & 1.6 & 1.8 & 0.8 & 0.1 & 0.9 & 0.8 & 0.9 \\
	Ma2015& 0.0 & 0.0 & 0.0 & 0.0 & 1.1 & 0.0 & 0.0 & 0.0 & 0.0 & 0.1 & 0.0 & 0.0 & 0.0 & 0.0 & 0.0 \\
	Mu2018& 0.0 & 0.0 & 0.0 & 0.0 & 0.2 & 0.0 & 0.0 & 0.0 & 0.0 & 0.0 & 0.0 & 0.0 & 0.0 & 0.0 & 0.0 \\
	\addlinespace[.2cm]
	\multicolumn{16}{l}{\quad\textit{Panel G: Sensitivity for $\mu_3$ }} \\
	Ours (full)& 81.1 & 42.3 & 78.0 & 58.4 & 59.5 & 83.4 & 45.0 & 83.7 & 56.6 & 54.6 & 78.2 & 33.8 & 73.2 & 40.8 & 44.0 \\
	Ours (partial)& 88.9 & 63.8 & 85.2 & 69.4 & 69.2 & 91.6 & 65.0 & 89.5 & 63.2 & 63.6 & 90.6 & 60.7 & 86.5 & 51.0 & 55.4 \\
	Ca2018& 78.3 & 37.8 & 82.1 & 74.9 & 79.5 & 83.6 & 30.2 & 85.4 & 81.1 & 84.2 & 86.2 & 18.5 & 87.2 & 84.3 & 87.2 \\
	Ho2018& 100.0 & 99.5 & 100.0 & 100.0 & 99.9 & 100.0 & 99.9 & 100.0 & 100.0 & 94.1 & 100.0 & 99.9 & 99.5 & 100.0 & 100.0 \\
	We2024& 3.3 & 0.5 & 3.7 & 3.0 & 3.3 & 1.6 & 0.2 & 1.9 & 1.6 & 1.8 & 0.8 & 0.1 & 0.9 & 0.8 & 0.9 \\
	Ma2015& 0.0 & 0.0 & 0.0 & 0.0 & 1.5 & 0.0 & 0.0 & 0.0 & 0.0 & 0.1 & 0.0 & 0.0 & 0.0 & 0.0 & 0.0 \\
	Mu2018& 0.0 & 0.0 & 0.0 & 0.0 & 0.5 & 0.0 & 0.0 & 0.0 & 0.0 & 0.0 & 0.0 & 0.0 & 0.0 & 0.0 & 0.0 \\
	\addlinespace[.2cm]
	\multicolumn{16}{l}{\quad\textit{Panel H: Sensitivity for $\mu_0$ }} \\
	Ours (full)& 82.1 & 42.4 & 79.0 & 56.5 & 60.3 & 83.2 & 45.5 & 83.0 & 55.8 & 56.5 & 76.3 & 34.0 & 78.1 & 40.8 & 44.5 \\
	Ours (partial)& 92.2 & 62.8 & 84.3 & 72.4 & 68.8 & 93.0 & 62.4 & 86.1 & 64.7 & 67.4 & 88.7 & 61.6 & 85.4 & 47.0 & 60.3 \\
	Ca2018& 83.7 & 43.6 & 84.8 & 80.5 & 80.8 & 86.4 & 34.6 & 87.0 & 83.2 & 85.0 & 87.3 & 26.8 & 88.0 & 85.2 & 88.0 \\
	Ho2018& 100.0 & 99.6 & 100.0 & 100.0 & 99.7 & 100.0 & 100.0 & 100.0 & 100.0 & 92.2 & 100.0 & 100.0 & 98.5 & 100.0 & 100.0 \\
	We2024& 3.2 & 0.6 & 3.6 & 2.9 & 3.4 & 1.6 & 0.2 & 1.8 & 1.5 & 1.8 & 0.8 & 0.1 & 0.9 & 0.8 & 0.9 \\
	Ma2015& 0.0 & 0.0 & 0.0 & 0.0 & 0.9 & 0.0 & 0.0 & 0.0 & 0.0 & 0.0 & 0.0 & 0.0 & 0.0 & 0.0 & 0.0 \\
	Mu2018& 0.0 & 0.0 & 0.0 & 0.0 & 0.3 & 0.0 & 0.0 & 0.0 & 0.0 & 0.0 & 0.0 & 0.0 & 0.0 & 0.0 & 0.0 \\
	\hline \hline 
 \end{tabular}
\end{footnotesize}
\end{table}
\begin{table}
\caption{Empirical specificity and sensitivity of various methods for MA errors, under different mean functions and error distributions.}
\begin{footnotesize}
\begin{tabular}{l| rrrrr | rrrrr | rrrrr}
	\hline \hline
	$n$ & \multicolumn{5}{c|}{50} & \multicolumn{5}{c|}{100} & \multicolumn{5}{c}{200} \\
	 $\mu$& $\Nc$ & $\Uc$ & $Exp$ & $Par_1$ & $Par_2$& $\Nc$ & $\Uc$ & $Exp$ & $Par_1$ & $Par_2$& $\Nc$ & $\Uc$ & $Exp$ & $Par_1$ & $Par_2$\\
	\hline
	\addlinespace[.2cm]
	\multicolumn{16}{l}{\quad\textit{Panel A: Specificity for $\mu_1$ }} \\
	Ours (full)& 86.9 & 93.9 & 86.4 & 90.7 & 88.9 & 87.2 & 96.1 & 86.4 & 92.0 & 90.6 & 86.8 & 97.0 & 85.9 & 93.9 & 92.8 \\
	Ours (partial)& 59.6 & 72.2 & 64.5 & 63.1 & 68.0 & 83.1 & 87.2 & 83.6 & 87.3 & 86.8 & 95.0 & 96.1 & 93.5 & 96.5 & 96.3 \\
	Ca2018& 99.9 & 100.0 & 99.9 & 99.7 & 99.7 & 100.0 & 100.0 & 99.9 & 99.8 & 99.8 & 100.0 & 100.0 & 99.9 & 99.9 & 99.9 \\
	Ho2018& 72.4 & 68.7 & 71.5 & 60.5 & 58.2 & 68.2 & 53.4 & 68.2 & 51.1 & 61.3 & 71.3 & 62.3 & 72.1 & 64.8 & 69.2 \\
	We2024& 99.7 & 99.7 & 99.7 & 99.7 & 99.7 & 99.8 & 99.8 & 99.9 & 99.8 & 99.8 & 99.9 & 99.9 & 99.9 & 99.9 & 99.9 \\
	Ma2015& 100.0 & 100.0 & 100.0 & 100.0 & 100.0 & 100.0 & 100.0 & 100.0 & 100.0 & 100.0 & 100.0 & 100.0 & 100.0 & 100.0 & 100.0 \\
	Mu2018& 100.0 & 100.0 & 100.0 & 100.0 & 100.0 & 100.0 & 100.0 & 100.0 & 100.0 & 100.0 & 100.0 & 100.0 & 100.0 & 100.0 & 100.0 \\
	\addlinespace[.2cm]
	\multicolumn{16}{l}{\quad\textit{Panel B: Specificity for $\mu_2$ }} \\
	Ours (full)& 86.7 & 94.7 & 86.7 & 91.4 & 88.9 & 87.1 & 96.6 & 86.5 & 93.3 & 90.9 & 86.8 & 97.0 & 86.1 & 94.6 & 92.9 \\
	Ours (partial)& 75.4 & 80.7 & 77.1 & 80.7 & 79.4 & 88.3 & 90.2 & 87.8 & 91.5 & 90.0 & 95.6 & 96.2 & 96.1 & 97.5 & 97.2 \\
	Ca2018& 99.9 & 100.0 & 99.9 & 99.7 & 99.7 & 100.0 & 100.0 & 99.9 & 99.8 & 99.8 & 100.0 & 100.0 & 99.9 & 99.9 & 99.9 \\
	Ho2018& 73.1 & 69.4 & 72.4 & 59.7 & 54.6 & 68.9 & 54.9 & 69.3 & 49.3 & 59.9 & 71.9 & 62.9 & 73.1 & 65.3 & 68.7 \\
	We2024& 99.7 & 99.7 & 99.7 & 99.7 & 99.7 & 99.8 & 99.8 & 99.9 & 99.8 & 99.8 & 99.9 & 99.9 & 99.9 & 99.9 & 99.9 \\
	Ma2015& 100.0 & 100.0 & 100.0 & 100.0 & 100.0 & 100.0 & 100.0 & 100.0 & 100.0 & 100.0 & 100.0 & 100.0 & 100.0 & 100.0 & 100.0 \\
	Mu2018& 100.0 & 100.0 & 100.0 & 100.0 & 100.0 & 100.0 & 100.0 & 100.0 & 100.0 & 100.0 & 100.0 & 100.0 & 100.0 & 100.0 & 100.0 \\
	\addlinespace[.2cm]
	\multicolumn{16}{l}{\quad\textit{Panel C: Specificity for $\mu_3$ }} \\
	Ours (full)& 86.1 & 93.5 & 86.5 & 91.1 & 88.6 & 86.8 & 96.1 & 86.6 & 93.1 & 90.9 & 86.7 & 97.1 & 86.1 & 94.6 & 92.9 \\
	Ours (partial)& 74.7 & 77.2 & 75.2 & 79.3 & 79.0 & 89.4 & 91.1 & 87.4 & 91.5 & 88.9 & 95.5 & 95.9 & 95.1 & 97.8 & 97.0 \\
	Ca2018& 99.9 & 99.8 & 99.8 & 99.7 & 99.7 & 99.9 & 99.9 & 99.9 & 99.8 & 99.8 & 100.0 & 100.0 & 99.9 & 99.9 & 99.9 \\
	Ho2018& 68.5 & 66.1 & 65.4 & 48.4 & 55.0 & 64.3 & 48.9 & 62.6 & 35.8 & 58.1 & 69.3 & 60.0 & 69.3 & 57.3 & 66.4 \\
	We2024& 99.7 & 99.7 & 99.7 & 99.7 & 99.7 & 99.8 & 99.9 & 99.9 & 99.8 & 99.8 & 99.9 & 99.9 & 99.9 & 99.9 & 99.9 \\
	Ma2015& 100.0 & 100.0 & 100.0 & 100.0 & 100.0 & 100.0 & 100.0 & 100.0 & 100.0 & 100.0 & 100.0 & 100.0 & 100.0 & 100.0 & 100.0 \\
	Mu2018& 100.0 & 100.0 & 100.0 & 100.0 & 100.0 & 100.0 & 100.0 & 100.0 & 100.0 & 100.0 & 100.0 & 100.0 & 100.0 & 100.0 & 100.0 \\
	\addlinespace[.2cm]
	\multicolumn{16}{l}{\quad\textit{Panel D: Specificity for $\mu_0$ }} \\
	Ours (full)& 87.0 & 94.5 & 86.6 & 91.3 & 89.4 & 87.1 & 96.2 & 86.9 & 92.8 & 91.3 & 86.5 & 97.1 & 86.4 & 94.5 & 92.6 \\
	Ours (partial)& 73.7 & 78.8 & 77.5 & 79.2 & 83.1 & 87.0 & 92.1 & 88.8 & 91.9 & 89.4 & 95.3 & 96.8 & 95.9 & 97.9 & 95.5 \\
	Ca2018& 99.9 & 100.0 & 99.9 & 99.8 & 99.7 & 100.0 & 100.0 & 99.9 & 99.8 & 99.8 & 100.0 & 100.0 & 99.9 & 99.9 & 99.9 \\
	Ho2018& 73.1 & 69.5 & 72.9 & 58.9 & 56.3 & 69.0 & 54.8 & 69.9 & 50.9 & 58.5 & 71.9 & 62.9 & 73.7 & 65.7 & 69.1 \\
	We2024& 99.7 & 99.7 & 99.7 & 99.7 & 99.7 & 99.8 & 99.8 & 99.9 & 99.8 & 99.8 & 99.9 & 99.9 & 99.9 & 99.9 & 99.9 \\
	Ma2015& 100.0 & 100.0 & 100.0 & 100.0 & 100.0 & 100.0 & 100.0 & 100.0 & 100.0 & 100.0 & 100.0 & 100.0 & 100.0 & 100.0 & 100.0 \\
	Mu2018& 100.0 & 100.0 & 100.0 & 100.0 & 100.0 & 100.0 & 100.0 & 100.0 & 100.0 & 100.0 & 100.0 & 100.0 & 100.0 & 100.0 & 100.0 \\
	\addlinespace[.2cm]
	\multicolumn{16}{l}{\quad\textit{Panel E: Sensitivity for $\mu_1$ }} \\
	Ours (full)& 73.2 & 29.7 & 70.7 & 46.6 & 56.3 & 85.4 & 37.6 & 83.0 & 57.5 & 62.0 & 92.8 & 50.1 & 90.9 & 62.5 & 59.4 \\
	Ours (partial)& 81.4 & 42.6 & 76.2 & 58.1 & 65.3 & 89.4 & 42.4 & 87.1 & 59.2 & 65.6 & 95.7 & 50.6 & 94.4 & 64.2 & 61.5 \\
	Ca2018& 82.2 & 39.6 & 83.5 & 78.9 & 80.1 & 86.0 & 34.2 & 86.1 & 82.8 & 84.3 & 86.9 & 24.3 & 87.3 & 84.6 & 87.6 \\
	Ho2018& 100.0 & 99.4 & 100.0 & 100.0 & 100.0 & 100.0 & 100.0 & 100.0 & 100.0 & 95.6 & 100.0 & 100.0 & 99.0 & 100.0 & 100.0 \\
	We2024& 2.9 & 0.5 & 3.6 & 3.0 & 3.4 & 1.5 & 0.3 & 1.8 & 1.5 & 1.7 & 0.7 & 0.1 & 0.9 & 0.8 & 0.9 \\
	Ma2015& 0.0 & 0.0 & 0.0 & 0.0 & 0.0 & 0.0 & 0.0 & 0.0 & 0.0 & 0.0 & 0.0 & 0.0 & 0.0 & 0.0 & 0.0 \\
	Mu2018& 0.0 & 0.0 & 0.0 & 0.0 & 0.2 & 0.0 & 0.0 & 0.0 & 0.0 & 0.0 & 0.0 & 0.0 & 0.0 & 0.0 & 0.0 \\
	\addlinespace[.2cm]
	\multicolumn{16}{l}{\quad\textit{Panel F: Sensitivity for $\mu_2$ }} \\
	Ours (full)& 83.7 & 40.3 & 82.4 & 59.9 & 64.0 & 90.1 & 44.9 & 88.8 & 62.2 & 62.2 & 94.4 & 51.9 & 93.2 & 60.0 & 61.5 \\
	Ours (partial)& 88.5 & 47.9 & 86.0 & 65.0 & 69.1 & 93.6 & 48.3 & 92.1 & 64.3 & 65.3 & 97.1 & 52.4 & 96.2 & 61.6 & 63.2 \\
	Ca2018& 83.7 & 42.6 & 83.8 & 80.3 & 80.9 & 86.2 & 34.8 & 86.8 & 83.2 & 85.0 & 87.1 & 25.5 & 87.8 & 85.0 & 88.0 \\
	Ho2018& 100.0 & 99.5 & 100.0 & 100.0 & 100.0 & 100.0 & 100.0 & 100.0 & 100.0 & 91.3 & 100.0 & 100.0 & 98.7 & 100.0 & 100.0 \\
	We2024& 2.9 & 0.5 & 3.6 & 2.9 & 3.4 & 1.5 & 0.2 & 1.8 & 1.5 & 1.8 & 0.8 & 0.1 & 0.9 & 0.8 & 0.9 \\
	Ma2015& 0.0 & 0.0 & 0.0 & 0.0 & 1.4 & 0.0 & 0.0 & 0.0 & 0.0 & 0.0 & 0.0 & 0.0 & 0.0 & 0.0 & 0.0 \\
	Mu2018& 0.0 & 0.0 & 0.0 & 0.0 & 0.5 & 0.0 & 0.0 & 0.0 & 0.0 & 0.0 & 0.0 & 0.0 & 0.0 & 0.0 & 0.0 \\
	\addlinespace[.2cm]
	\multicolumn{16}{l}{\quad\textit{Panel G: Sensitivity for $\mu_3$ }} \\
	Ours (full)& 84.1 & 42.4 & 81.7 & 60.9 & 65.7 & 90.7 & 46.3 & 87.6 & 62.1 & 62.5 & 94.2 & 50.9 & 92.3 & 61.0 & 61.8 \\
	Ours (partial)& 88.3 & 50.6 & 86.5 & 66.5 & 69.8 & 94.2 & 48.6 & 91.1 & 64.4 & 65.7 & 97.0 & 51.5 & 95.7 & 62.5 & 63.6 \\
	Ca2018& 77.6 & 33.8 & 82.2 & 74.9 & 79.8 & 84.1 & 25.9 & 85.5 & 81.0 & 84.6 & 85.9 & 18.7 & 87.1 & 84.1 & 87.6 \\
	Ho2018& 100.0 & 99.2 & 100.0 & 100.0 & 99.8 & 100.0 & 99.9 & 100.0 & 100.0 & 94.7 & 100.0 & 99.9 & 99.4 & 100.0 & 100.0 \\
	We2024& 3.0 & 0.5 & 3.5 & 2.8 & 3.4 & 1.5 & 0.2 & 1.8 & 1.5 & 1.8 & 0.7 & 0.1 & 0.9 & 0.8 & 0.9 \\
	Ma2015& 0.0 & 0.0 & 0.0 & 0.0 & 1.4 & 0.0 & 0.0 & 0.0 & 0.0 & 0.1 & 0.0 & 0.0 & 0.0 & 0.0 & 0.0 \\
	Mu2018& 0.0 & 0.0 & 0.0 & 0.0 & 0.5 & 0.0 & 0.0 & 0.0 & 0.0 & 0.0 & 0.0 & 0.0 & 0.0 & 0.0 & 0.0 \\
	\addlinespace[.2cm]
	\multicolumn{16}{l}{\quad\textit{Panel H: Sensitivity for $\mu_0$ }} \\
	Ours (full)& 83.9 & 41.6 & 82.7 & 59.7 & 62.5 & 90.2 & 47.4 & 88.5 & 63.1 & 64.3 & 94.1 & 53.2 & 92.8 & 62.9 & 61.2 \\
	Ours (partial)& 91.5 & 50.3 & 86.1 & 66.9 & 66.5 & 94.5 & 48.8 & 90.9 & 65.9 & 67.6 & 97.2 & 53.1 & 95.3 & 64.2 & 64.2 \\
	Ca2018& 83.7 & 42.0 & 84.5 & 79.9 & 81.2 & 86.3 & 36.0 & 86.9 & 83.0 & 85.5 & 87.3 & 26.6 & 87.9 & 85.3 & 88.0 \\
	Ho2018& 100.0 & 99.4 & 100.0 & 100.0 & 99.6 & 100.0 & 99.9 & 100.0 & 100.0 & 92.3 & 100.0 & 100.0 & 98.4 & 100.0 & 100.0 \\
	We2024& 3.0 & 0.5 & 3.6 & 2.9 & 3.3 & 1.5 & 0.2 & 1.8 & 1.5 & 1.8 & 0.7 & 0.1 & 0.9 & 0.8 & 0.9 \\
	Ma2015& 0.0 & 0.0 & 0.0 & 0.0 & 0.6 & 0.0 & 0.0 & 0.0 & 0.0 & 0.0 & 0.0 & 0.0 & 0.0 & 0.0 & 0.0 \\
	Mu2018& 0.0 & 0.0 & 0.0 & 0.0 & 0.3 & 0.0 & 0.0 & 0.0 & 0.0 & 0.0 & 0.0 & 0.0 & 0.0 & 0.0 & 0.0 \\
	\hline \hline 
 \end{tabular}
\end{footnotesize}
\end{table}